\documentclass[11pt]{amsart}

\usepackage{amssymb}
\usepackage{epsf}
\usepackage{a4wide}
\usepackage[normal]{caption2}
\usepackage{dsfont}

\makeatletter
\newcommand\tabcaption{\def\@captype{table}\caption}
\makeatother

\newtheorem{theorem}{Theorem}
\newtheorem{prop}{Proposition}
\newtheorem{lem}{Lemma}
\newtheorem{coro}{Corollary}

\newcommand{\ABC}{\mathcal{A}}
\newcommand{\bs}[1]{\boldsymbol{#1}}
\newcommand{\openset}[1]{\overset{\circ}{#1}}
\newcommand{\diamondset}[1]{\overset{\diamond}{#1}}
\newcommand{\NN}{\mathbb{N}}
\newcommand{\ZZ}{\mathbb{Z}}
\newcommand{\QQ}{\mathbb{Q}}
\newcommand{\RR}{\mathbb{R}}
\newcommand{\CC}{\mathbb{C}}
\newcommand{\RE}{\operatorname{Re}}
\newcommand{\IM}{\operatorname{Im}}
\newcommand{\CPS}{\Sigma\textnormal{Kol}(3,1)}

\begin{document}

\title{Kolakoski-(3,1) is a (deformed) model set}  

\author{Michael Baake}
\author{Bernd Sing}
\address{Institut f\"{u}r Mathematik, Universit\"{a}t Greifswald,
Jahnstr.~15a, 17487 Greifswald, Germany}
\email{mbaake@uni-greifswald.de}
\email{sing@uni-greifswald.de}
\urladdr{ http://schubert.math-inf.uni-greifswald.de }

\begin{abstract} 
Unlike the (classical) Kolakoski sequence on the alphabet $\{1,2\}$, its
analogue on $\{1,3\}$ can be related to a primitive substitution rule. Using
this connection, we prove that the corresponding bi-infinite fixed point is a
regular generic model set and thus has a pure point diffraction spectrum. The
Kolakoski-$(3,1)$ sequence is then obtained as a deformation, without losing
the pure point diffraction property.
\end{abstract}
\maketitle

\section{Introduction}
 
A one-sided infinite sequence $\omega$ over the alphabet $\ABC=\{1,2\}$ is
called a (classical) \textit{Kolakoski sequence} (named after W.~Kolakoski who
introduced it in 1965, see~\cite{Kol65}), if it equals the sequence defined by
its run lengths, e.g.:
\begin{equation}\label{eq:kol}
\begin{array}{ccccccccccccccc}
\omega & = & \underbrace{22} & \underbrace{11} & \underbrace{2} &
\underbrace{1} & \underbrace{22} & \underbrace{1} & \underbrace{22} &
\underbrace{11} & \underbrace{2} & \underbrace{11} & \ldots && \\
&& 2 & 2 & 1 & 1 & 2 & 1 & 2 & 2 & 1 & 2 & \ldots & = & \omega.
\end{array}
\end{equation}
Here, a \textit{run} is a maximal subword consisting of identical letters. The
sequence $\omega'=1\omega$ is the only other sequence which has this property.

One way to obtain $\omega$ of (\ref{eq:kol}) is by starting with $2$ as a seed
and iterating the two substitutions 
\begin{equation*}
\sigma^{}_0: \begin{array}{lcl} 1 & \mapsto & 2 \\ 2 & \mapsto & 22 \end{array}
\quad \text{and} \quad
\sigma^{}_1: \begin{array}{lcl} 1 & \mapsto & 1 \\ 2 & \mapsto & 11,
\end{array} 
\end{equation*}
alternatingly, i.e., $\sigma^{}_0$ substitutes letters on even positions and
$\sigma^{}_1$ letters on odd positions (we begin counting at $0$):
\begin{equation*}
2 \mapsto 22 \mapsto 2211 \mapsto 221121 \mapsto 221121221 \mapsto \ldots
\end{equation*}
Clearly, the iterates converge to the Kolakoski sequence $\omega$ (in the
obvious product topology), and $\omega$ is the unique (one-sided) fixed point
of this iteration.  

One can generalize this by choosing a different alphabet
$\ABC=\{p,q\}$ (we are only looking at alphabets with $\operatorname{card}
(\ABC) = 2$), e.g., $\ABC=\{1,3\}$, which is the main focus of this paper.
Such a (generalized) Kolakoski sequence, which is also equal to the sequence of
its run lengths, can be obtained by iterating the two
substitutions 
\begin{equation*}
\sigma^{}_0: \begin{array}{lcl} q & \mapsto & p^q \\ p & \mapsto & p^p
\end{array} \quad \text{and} \quad
\sigma^{}_1: \begin{array}{lcl} q & \mapsto & q^q \\ p & \mapsto & q^p
\end{array} 
\end{equation*}
alternatingly. Here, the starting letter of the sequence is $p$. We will call
such a sequence Kolakoski-$(p,q)$ sequence, or Kol$(p,q)$ for short. The
classical Kolakoski sequence $\omega$ of (\ref{eq:kol}) is therefore denoted
by Kol$(2,1)$ (and $\omega'$ by Kol$(1,2)$).  

While little is known about the classical Kolakoski sequence
(see~\cite{Dek97}), and the same holds for all Kol$(p,q)$ with $p$ odd and $q$
even or vice versa (see~\cite{Diplom}), the situation is more favourable if
$p$ and $q$ are either both even or both odd. If both are even, one can rewrite
the substitution as a substitution of constant length by building blocks of 4
letters (see~\cite{Diplom, InPrep}). Spectral properties can then be deduced
by a criterion of Dekking~\cite{Dek78}. The case where both symbols are
odd will be studied in this paper exemplarily on Kol$(3,1)$. 

It is our aim to determine structure and order
of the sequence Kol$(3,1)$. This will require two steps: First, we relate it
to a unimodular substitution of Pisot type and prove that the corresponding
aperiodic point set is a regular generic model set. Second, we relate this
back to the original Kol$(3,1)$ by a deformation. Here, the first step is a
concrete example of the general conjecture that all unimodular substitutions
of Pisot type are regular model sets (however, not always generic). This
general conjecture cannot be proved by an immediate application of our
strategy, but we hope that our method sheds new light on it.

\smallskip
\noindent
\textsc{Remark}: Every Kol$(p,q)$ can uniquely be extended to a
bi-infinite (or two-sided) sequence. The one-sided sequence (to the right) is
Kol$(p,q)$ as explained above. The added part to the left is a reversed copy
of Kol$(q,p)$, e.g., in the case of the classical Kolakoski sequence of
(\ref{eq:kol}), this reads as
\begin{equation*}
\ldots 11221221211221|22112122122112 \ldots,
\end{equation*}
where ``$|$'' denotes the seamline between the one-sided sequences. Note that,
if $q=1$ (or $p=1$), the bi-infinite sequence is mirror symmetric around the
first position to the left (right) of the seamline. The bi-infinite sequence
equals the sequence of its run lengths, if counting is begun at the seamline. 
Alternatively, one can get such a bi-infinite sequence by starting with $q|p$
and applying the two substitutions to get $\sigma^{}_1(q)|\sigma^{}_0(p)$ in
the first step and so forth. This also implies that Kol$(p,q)$ and Kol$(q,p)$
will have the same spectral properties, and it suffices to study one of them.

\section{Kol$(3,1)$ as substitution}

If both letters are odd numbers, one can build blocks of $2$ letters and obtain
an (ordinary) substitution. Setting\footnote
{
That Kol$(3,1)$ can be related to a substitution is well-known, e.g.,
in~\cite{Dek80}, a substitution over an alphabet with four letters is given,
while~\cite{Sir98} uses the same substitution with three letters as we do. We
thank the referee for pointing this last reference out to us. 
} 
$A=33$, $B=31$ and $C=11$ in the case of Kol$(3,1)$, this substitution
$\sigma$ and its \textit{substitution matrix} $\bs{M}$ (sometimes called
\textit{incidence matrix} of the substitution) are given by
\begin{equation} \label{eq:subs}
\sigma: \begin{array}{lcl} A & \mapsto & ABC \\ B & \mapsto & AB \\ C &
  \mapsto & B \end{array} \quad \text{and} \quad \bs{M}=\left(
  \begin{array}{ccc} 1 & 1 & 1 \\ 1 & 1 & 0 \\ 0 & 1 & 0 \end{array} \right),
\end{equation}
where the entry $M_{ij}$ is the number of occurrences of $j$ in $\sigma(i)$
($i,j \in \{A,B,C\}$; sometimes the transposed matrix is used). A bi-infinite
fixed point can be obtained as follows:
\begin{equation}\label{eq:BA}
B|A \, \mapsto \, AB|ABC \, \mapsto \, ABCAB|ABCABB \, \mapsto \, \ldots
\end{equation}
This corresponds to 
\begin{equation}\label{eq:31}
\ldots3331113331|333111333131\ldots
\end{equation}
which is the unique bi-infinite Kol$(3,1)$ according to our above convention. 
The matrix $\bs{M}$ is primitive because $\bs{M}^3$ has positive entries
only. The characteristic polynomial $P(x)$ of $\bs{M}$ is 
\begin{equation}\label{eq:charpol}
P(x) := \det(x\,\mathds{1} -\bs{M}) = x^3 - 2 \, x^2 - 1,
\end{equation}  
which is irreducible over $\ZZ$ (there is no solution $\bmod \, 3$) and over
$\QQ$ (every rational algebraic integer is an integer). The discriminant $D$ of
$P(x)$ is $D = \frac{59}{108}$, so $P(x)$ has one real root
$\alpha$ and two complex conjugate roots $\beta$ and $\overline{\beta}$. One
gets 
\begin{equation*}
2.21 \approx \alpha > 1 > |\beta| \approx 0.67 >0,
\end{equation*}
wherefore $\alpha$ is a \textit{Pisot-Vijayaraghavan number} (i.e., an
algebraic integer greater than $1$ whose algebraic conjugates are all 
less than $1$ in modulus), and $\sigma$ is a substitution of \textit{Pisot
  type}. Since $\det(\bs{M}) = 1$, the roots $\alpha$, $\beta$ and
$\overline{\beta}$ are also \textit{algebraic units}, and the associated
substitution is said to be \textit{unimodular}. Note that $\RE(\beta) =
1-\frac{\alpha}{2}$. If necessary, we will choose $\beta$ such that
$\IM(\beta)>0$ in the following calculations (the other possibility only leads
to overall minus signs). 

There is a natural geometric representation of such a substitution by
inflation, compare~\cite{LGJJ93}. Here, one associates bond lengths (or 
intervals) $\ell^{}_{A}$, $\ell^{}_{B}$ and $\ell^{}_{C}$ to each
letter. These bond lengths are given by the components of the right eigenvector
which belongs to the (real) eigenvalue $\alpha$ and is unique (up to
normalization) by the Perron-Frobenius theorem. The normalization can be
chosen so that
\begin{equation*}
\ell^{}_{A} = \alpha^2 - \alpha \approx 2.66, \qquad \ell^{}_{B} = \alpha
\approx 2.21 \quad \text{and} \quad \ell^{}_{C} = 1.
\end{equation*}
Inflating  the bond lengths by a factor of $\alpha$ and dividing them into
original intervals just corresponds to the substitution (because $\alpha \cdot
\ell^{}_{A} = \ell^{}_{A} + \ell^{}_{B} + \ell^{}_{C}$, etc.). We will denote
this realization of the bi-infinite fixed point with natural bond lengths 
(respectively the point set associated with this realization where we mark the
left endpoints of the intervals by their name) by $\CPS$, reserving
``Kol$(3,1)$'' for the case of unit (or integer) bond lengths. 

On the other hand, the frequencies $\rho^{}_{A}$, $\rho^{}_{B}$ and
$\rho^{}_{C}$ of the letters in the infinite sequence are given by the
components of the left eigenvector of $\bs{M}$ to the eigenvalue
$\alpha$. This gives 
\begin{equation*}
\rho^{}_{A} = \frac12\,(-\alpha^2+3\,\alpha-1) \approx 0.38, \quad
\rho^{}_{B} = \alpha^2-2\,\alpha \approx 0.45, \quad 
\rho^{}_{C} = \frac12\,(-\alpha^2+\alpha+3) \approx 0.17,
\end{equation*}
with $\rho^{}_{A} + \rho^{}_{B} + \rho^{}_{C} = 1$.
Therefore, the average bond length $\ell$ in the geometric representation is
\begin{equation}\label{eq:avlen}
\ell = \rho^{}_{A} \cdot \ell^{}_{A} + \rho^{}_{B} \cdot \ell^{}_{B} +
\rho^{}_{C} \cdot \ell^{}_{C} = \frac12\,(-\alpha^2+\alpha+7)\approx 2.17, 
\end{equation}
and the frequencies of $3$s and $1$s in Kol$(3,1)$ can easily be calculated
to be $\rho^{}_{3} = \frac12 (\alpha-1) \approx 0.60$ and $\rho^{}_{1} =
\frac12(-\alpha+3) \approx 0.40$.  

\smallskip
\noindent
\textsc{Remark}: In the case where $p$ and $q$ are odd (positive) integers,
one gets unimodular substitutions of Pisot type iff $p = q \pm 2$. More
generally, one gets substitutions of Pisot type iff $2 \cdot (p+q) \ge
(p-q)^2$ holds. Otherwise, all the eigenvalues are greater than $1$ in
modulus, see~\cite{Diplom}. 

\section{Model Set and IFS \label{sec:IFS}}

A \textit{model set} $\varLambda(\Omega)$ (or \textit{cut-and-project set})
in \textit{physical space} $\RR^d$ is defined within the following general
cut-and-project scheme~\cite{Moo00,Baa02}
\begin{equation} \label{eq:cutproj}
\renewcommand{\arraystretch}{1.5}
\begin{array}{ccccc}
& \pi & & \pi_{\textnormal{int}}^{} & \vspace*{-1.5ex} \\
\RR^{d} & \longleftarrow & \RR^{d}\times H & \longrightarrow & H \\
 & \mbox{\raisebox{-1.5ex}{\footnotesize
     \textnormal{1--1}}}\!\!\!\!\nwarrow\;\; & \cup &  
\;\;\nearrow\!\!\!\!\mbox{\raisebox{-1.5ex}{\footnotesize \textnormal{dense}}}
& \\ &  & \varGamma &  & 
\end{array}
\end{equation}
where the \textit{internal space} $H$ is a locally compact Abelian group,
and $\varGamma\subset\RR^{d}\times H$ is a \textit{lattice}, i.e., a co-compact
discrete subgroup of $\RR^{d}\times H$. The projection
$\pi_{\textnormal{int}}^{}(\varGamma)$ is assumed to be dense in internal
space, and the projection $\pi$ into physical space has to be one-to-one on
$\varGamma$. The model set $\varLambda(\Omega)$ is 
\begin{equation*}
\varLambda(\Omega)\;=\;\left\{\pi(x)\mid x\in\varGamma,\,
\pi_{\textnormal{int}}^{}(x)\in\Omega \right\}\;\subset\; \RR^{d},
\end{equation*}
where the \textit{window} $\Omega\subset H$ is a relatively compact set with
non-empty interior. If we set $L = \pi(\varGamma) \subset \RR^{d}$, we can
define, for $x \in L$, the \textit{star  map} $^{\star}:\,L \mapsto H$ by
$x^{\star} = \pi_{\textnormal{int}}^{} \circ \left(\pi|^{}_{L}\right)^{-1}(x)$,
see~\cite{BM00}. So we have $\varGamma = \{(x, x^{\star}) \mathbin| x \in
L\}$ and $L^{\star} = \pi_{\textnormal{int}}^{}(\varGamma)$.  If the boundary
$\partial\Omega$ of the window has vanishing Haar measure in $H$, we call
$\varLambda(\Omega)$ a \textit{regular model set}. If, in addition,
$\pi_{\textnormal{int}}(\varGamma)\cap\partial\Omega = \varnothing$, the
model set is called non-singular or \textit{generic}.

Every model set is also a \textit{Delone set} (or \textit{Delaunay set}),
i.e., it is both \textit{uniformly discrete}\footnote
{
A set $\varLambda$ is \textit{uniformly discrete} if $\exists r>0$ s.t.\ every
open ball of radius $r$ contains at most one point of $\varLambda$. 
} 
and \textit{relatively dense}\footnote
{
A set $\varLambda$ is \textit{relatively dense} if $\exists R>0$ s.t.\ every
closed ball of radius $R$ contains at least one point of $\varLambda$.
}. 
A Delone set $X$ is a \textit{Meyer set}, if also $X - X$ is a Delone
set. Every model set is a Meyer set, see~\cite{Moo97}. 

We will now construct a model set $\varLambda(\Omega)$ and -- in a first
step -- show that this model set differs from $\CPS$ at most on
positions of density $0$. By Galois conjugation (see~\cite{LGJJ93, BS66}),
which here corresponds to the star map as we will see, we find a lattice 
\begin{equation}\label{eq:lattice}
\varGamma \, = \,
\ZZ\cdot\bs{v}^{}_{A}+\ZZ\cdot\bs{v}^{}_{B}+\ZZ\cdot\bs{v}^{}_{C}
\, \subset\, \RR\times\CC \, \simeq \, \RR^3
\end{equation}
where 
\begin{equation*}
\bs{v}^{}_{A} = \left( \begin{array}{c} \alpha^2-\alpha \\ \RE(\beta^2-\beta)
    \\ \IM(\beta^2-\beta) \end{array} \right), \qquad \bs{v}^{}_{B} = \left(
    \begin{array}{c} \alpha \\ \RE(\beta) \\ \IM(\beta) \end{array} \right)
    \quad \text{and} \quad \bs{v}^{}_{C} = \left( \begin{array}{c} 1 \\ 1 \\ 0
    \end{array} \right). 
\end{equation*}  
The projection $\pi$ (i.e., the projection on the first coordinate) is
injective on $\varGamma$ because $\QQ(\alpha)$ is a $\QQ$-vector space of
dimension $3$ with ($\QQ$-)linearly independent elements $1$, $\alpha$ and
$\alpha^2$. Also, $\pi(\varGamma)=\ZZ[\alpha]$ is dense.
To see that $\pi_{\textnormal{int}}^{}(\varGamma)$ is dense, we
note that $\pi_{\textnormal{int}}^{}(\varGamma) = \ZZ[\beta]$ and that $1$ and
$\beta$ are linearly independent. So, $\beta_{}^n$ and $\beta_{}^{n+1}$ are
also linearly independent for all $n \in \NN$, and their $\ZZ$-span forms a
two-dimensional lattice in $\CC$, which is a uniformly discrete subset of
$\ZZ[\beta]$. Since $|\beta|<1$, one can choose, for every $\varepsilon > 0$,
an $n$, such that there is a lattice point (of the lattice $\ZZ\cdot\beta_{}^n
+ \ZZ\cdot\beta_{}^{n+1}$) in every ball of radius $\varepsilon$, so
$\ZZ[\beta]$ is dense in $\CC$. Note, that $\pi_{\textnormal{int}}$ is also
injective on $\varGamma$ (this can be seen from $\RE(\beta) =
1-\frac{\alpha}{2}$ and $\RE(\beta^2) = 2 - \frac{\alpha^2}{2}$).
So we have established:

\begin{prop}\label{prop:cps_scheme}
With $\CC \simeq \RR^2$, $\varGamma$ of $(\ref{eq:lattice})$ and the natural
projections $\pi$ and $\pi_{\textnormal{int}}^{}$, we obtain the following
cut-and-project scheme: 
\begin{equation} \label{eq:cutproj1}
\renewcommand{\arraystretch}{1.5}
\begin{array}{ccccc}
& \pi & & \pi_{\textnormal{int}}^{} & \vspace*{-1.5ex} \\
\RR & \longleftarrow & \RR\times\RR^2  & \longrightarrow & \RR^2 \\
 & \mbox{\raisebox{-1.5ex}{\footnotesize 
    ${\begin{subarray}{c} \textnormal{dense} \\ \textnormal{1--1}
     \end{subarray}}$
}}\!\!\!\!\nwarrow\;\; & \cup &   
\;\;\nearrow\!\!\!\!\mbox{\raisebox{-1.5ex}{\footnotesize 
    ${\begin{subarray}{c} \textnormal{dense} \\ \textnormal{1--1}
    \end{subarray}}$}} & \\ &  & \varGamma &  & 
\end{array}
\end{equation}
Furthermore, we have
\begin{equation*}
\pi(\varGamma) = \ZZ[\alpha] \quad \text{and} \quad
\pi_{\textnormal{int}}^{}(\varGamma) = \ZZ[\beta], 
\end{equation*}
where $\alpha$ is the real root of $(\ref{eq:charpol})$ and $\beta$ one of the
complex conjugate ones.\qed
\end{prop}

In order to describe $\CPS$, the main task is now to determine the
appropriate windows $\Omega^{}_{A}$, $\Omega^{}_{B}$ and
$\Omega^{}_{C}$ (one for each letter; 
$\Omega=\Omega^{}_{A}\cup\Omega^{}_{B}\cup\Omega^{}_{C}$). For
these windows, the substitution rule $\sigma$ of (\ref{eq:subs}) induces the
following \textit{iterated function system} (IFS for short) in internal
space\footnote
{
For later reference, we write: 
\begin{equation*}\tag{\ref{eq:IFS}'}
\begin{array}{lclclcl}
\Omega^{}_{A} & = & f^{}_1(\Omega^{}_{A}) & \cup &
f^{}_1(\Omega^{}_{B}) \\ 
\Omega^{}_{B} & = & f^{}_3(\Omega^{}_{A}) & \cup &
f^{}_3(\Omega^{}_{B}) & \cup & f^{}_1(\Omega^{}_{C}) \\
\Omega^{}_{C} & = & f^{}_0(\Omega^{}_{A}), \\
\end{array}
\end{equation*}
where $f^{}_1$ and $f^{}_3$ are defined as in (\ref{eq:f-fcns}), and $f^{}_0
(z) = \beta \, z + \beta^2$.
}, 
cf.~\cite{LGJJ93}: 
\begin{equation} \label{eq:IFS} 
\begin{array}{lclclcl}
\Omega^{}_{A} & = & \beta \, \Omega^{}_{A} & \cup & \beta \,
\Omega^{}_{B} \\ 
\Omega^{}_{B} & = & \beta \, \Omega^{}_{A} + \beta^2 - \beta & \cup &
\beta \, \Omega^{}_{B} + \beta^2 - \beta & \cup & \beta \, \Omega^{}_{C}
\\ 
\Omega^{}_{C} & = & \beta \, \Omega^{}_{A} + \beta^2. \\
\end{array}
\end{equation}
This IFS is obtained as follows: We denote by $\varLambda^{}_{A}$ the subset of
$\CPS$ of left endpoints of intervals of type $A$ (of length $\ell^{}_{A}$),
and similar for $\varLambda^{}_{B}$ and $\varLambda^{}_{C}$ (we have $\CPS =
\varLambda^{}_{A} \dot{\cup} \varLambda^{}_{B} \dot{\cup} \varLambda^{}_{C}$,
where $\dot{\cup}$ denotes disjoint union). Then, the substitution $\sigma$
of~(\ref{eq:subs}) induces the following equations for these Delone sets in
$\RR$: 
\begin{equation}\label{eq:IFS_physical}
\begin{array}{lclclcl}
\varLambda^{}_{A} & = & \alpha \, \varLambda^{}_{A} & \cup & \alpha \,
\varLambda^{}_{B} \\ 
\varLambda^{}_{B} & = & \alpha \, \varLambda^{}_{A} + \alpha^2 - \alpha & \cup
& \alpha \, \varLambda^{}_{B} + \alpha^2 - \alpha & \cup & \alpha \,
\varLambda^{}_{C} \\ 
\varLambda^{}_{C} & = & \alpha \, \varLambda^{}_{A} + \alpha^2. \\
\end{array}
\end{equation}
Applying the star map to these equations yields~(\ref{eq:IFS}). In this sense,
the iteration of the IFS~(\ref{eq:IFS}) in internal space corresponds to the
iteration~(\ref{eq:IFS_physical}) in physical space (note that by
Proposition~\ref{prop:cps_scheme} the star map is bijective on $\ZZ[\alpha]$).

\begin{figure}[t]
\begin{minipage}[c]{0.5\textwidth}
\centerline{\epsfxsize=\textwidth\epsfbox{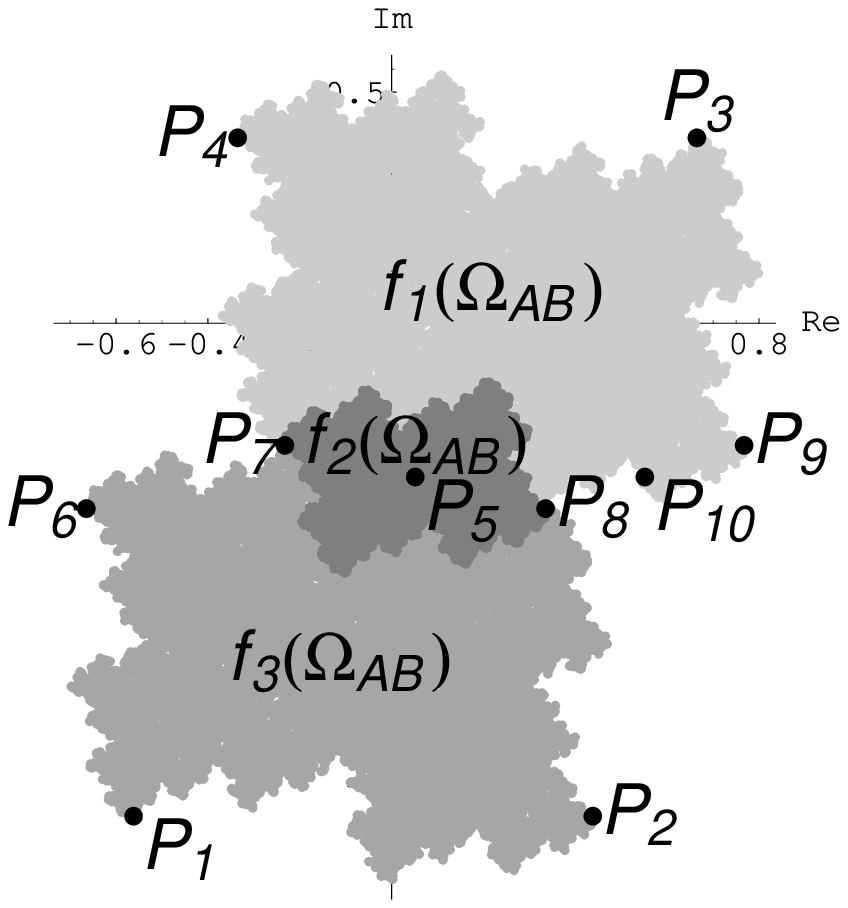}}
\setcaptionwidth{\textwidth}
\caption{Components $f^{}_{i}(\Omega^{}_{AB})$ of the set
  $\Omega^{}_{AB}$ and the points $P^{}_{1}, \ldots, P^{}_{10}$ used in the
  proof of Proposition~\ref{prop:jt}. $\Omega^{}_{AB}$ is the union of the
  three shaded areas.\label{fig:win1}} 
\end{minipage}\hfill
\begin{minipage}[c]{0.4\textwidth}
\centerline{\begin{tabular}{|c|c|}\hline
Point & Coordinate \\ \hline \hline \rule[-1.5ex]{0ex}{4ex}
$P^{}_1$ & $\frac12 \, (\beta^2 - 3 \, \beta - 1)$ \\ \rule[-1.5ex]{0ex}{4ex}
$P^{}_2$ & $\frac12 \, (\beta^2 - 3 \, \beta + 1)$ \\ \rule[-1.5ex]{0ex}{4ex}
$P^{}_3$ & $\frac12 \, (- \beta^2 + \beta + 1)$ \\ \rule[-1.5ex]{0ex}{4ex}
$P^{}_4$ & $\frac12 \, (- \beta^2 + \beta - 1)$ \\ \rule[-1.5ex]{0ex}{4ex}
$P^{}_5$ & $- \frac12 \, \beta$ \\ \rule[-1.5ex]{0ex}{4ex}
$P^{}_6$ & $\frac12 \, (\beta^2 - \beta - 1)$ \\ \rule[-1.5ex]{0ex}{4ex}
$P^{}_7$ & $\frac12 \, (- \beta^2 - \beta - 1)$ \\ \rule[-1.5ex]{0ex}{4ex}
$P^{}_8$ & $\frac12 \, (\beta^2 - \beta + 1)$ \\ \rule[-1.5ex]{0ex}{4ex}
$P^{}_9$ & $\frac12 \, (- \beta^2 - \beta + 1)$ \\ \rule[-1.5ex]{0ex}{4ex}
$P^{}_{10}$ & $\frac12 \, (- \beta + 1)$ \\ \hline 
\end{tabular}}
\vspace{1.5ex}
\tabcaption{Coordinates of the points $P^{}_{1}, \ldots, P^{}_{10}$ in
  Figure~\ref{fig:win1}. \label{tab:points}}
\end{minipage}
\end{figure}

\begin{figure}[t]
\centerline{\epsfxsize=0.5\textwidth\epsfbox{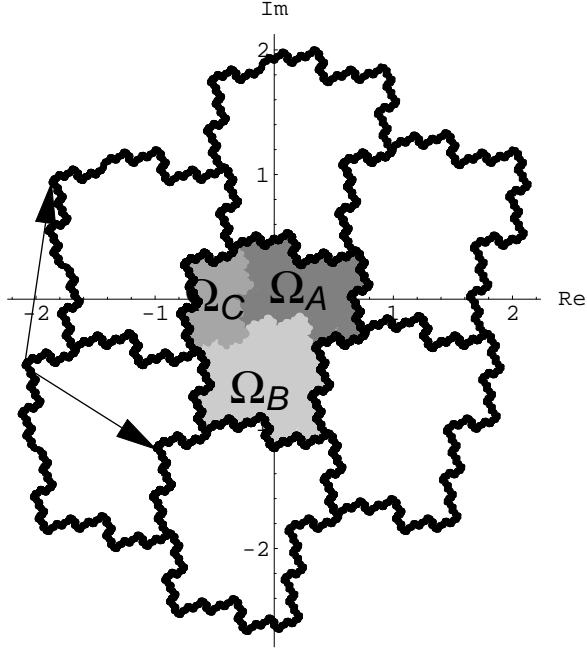}}
\caption{Components of the set $\Omega$, the periodic tiling of the plane
  ($\CC$) with it and the corresponding translation vectors. \label{fig:win2}} 
\end{figure}

Setting $\Omega^{}_{AB} = \Omega^{}_{A} \cup \Omega^{}_{B}$ in~(\ref{eq:IFS}),
the system decouples and we remain with the simpler IFS
\begin{equation} \label{eq:IFS-AB}
\begin{array}{lclclcl}
\Omega^{}_{AB} & = & f^{}_1(\Omega^{}_{AB}) & \cup &
f^{}_2(\Omega^{}_{AB}) & \cup & f^{}_3(\Omega^{}_{AB}) \\ 
\end{array}
\end{equation}
where 
\begin{equation} \label{eq:f-fcns}
f^{}_1(z) = \beta \, z, \qquad f^{}_2(z) = \beta^3 \, z + \beta^3 \quad
\text{and} \quad f^{}_3(z) = \beta \, z + \beta^2 - \beta.
\end{equation}
The mappings $f^{}_{i}: \CC \mapsto \CC$
are contractions ($|\beta|<1$), so that Hutchinson's theorem~\cite[Section
3.1(3)]{Hut81} guarantees a unique compact solution of (\ref{eq:IFS-AB}),
called the \textit{attractor} of the IFS. The sets $\Omega^{}_{A}$,
$\Omega^{}_{B}$ and $\Omega^{}_{C}$ can be calculated from
$\Omega^{}_{AB}$ as 
\begin{equation}\label{eq:affine}
\Omega^{}_{A} = f^{}_1(\Omega^{}_{AB}), \qquad \Omega^{}_{B} =
f^{}_2(\Omega^{}_{AB}) \cup f^{}_3(\Omega^{}_{AB}), \quad \text{and}
\quad \Omega^{}_{C} = f^{}_4(\Omega^{}_{AB}),
\end{equation}
where $f^{}_4(z) = \beta^2 \, z + \beta^2$. They are also compact sets in the
plane. For the components of $\Omega^{}_{AB}$, see Figure~\ref{fig:win1}; the
windows $\Omega^{}_{A}$, etc., are shown in Figure~\ref{fig:win2}. Note that
the decoupling of the IFS (i.e., the step from~(\ref{eq:IFS})
to~(\ref{eq:IFS-AB})) lies at the heart of our argument and seems to be the
reason that we cannot immediately generalize our method to other unimodular
substitutions of Pisot type\footnote
{
The substitution~(\ref{eq:subs}) can be analyzed by the balanced pair
algorithm as described in~\cite{SS02}. This algorithm also confirms that it
has pure point spectrum, but one does not get the model set property. 
}, 
because no such decoupling emerges in general. 

The \textit{similarity dimension} $s$ of a set given by an IFS is the unique
non-negative number $s$ such that the contraction constants to the power of $s$
add up to $1$ (see~\cite{Edg90}). For $\Omega^{}_{AB}$, this means
\begin{equation*}
|\beta|^{3s} + 2 \, |\beta|^{s} = 1
\end{equation*}
with solution $s=s(\Omega^{}_{AB})=2$ (because $\alpha \cdot |\beta|^2 = 1$,
the substitution is unimodular). The similarity dimension
$s$ of a set is connected to its \textit{Hausdorff dimension} $h$ by $h \le s$
where equality holds if the \textit{open set condition} (OSC for short) is
satisfied~\cite{Edg90}. An IFS with mappings $f^{}_{i}$ satisfies the OSC iff
there exists a nonempty open set $U$ such that $f^{}_{i}(U) \cap f^{}_{j}(U) =
\varnothing$  for $i \neq j$ and $f^{}_{i}(U) \subset U$ for all $i$. It is
easy to see that the corresponding self-similar set $\Omega^{}_{AB}$ must be
contained in the closure $\overline{U}$ so that the pieces
$f^{}_{i}(\Omega^{}_{AB}) \subset \overline{f^{}_{i}(U)}$ can intersect
at their boundaries but cannot have interior points in common~\cite{Ban97}. If 
their boundaries do intersect, the IFS is called \textit{just touching}. 

\begin{prop} \label{prop:jt} 
The IFS of $(\ref{eq:IFS-AB})$ for $\Omega^{}_{AB}$ is just touching. 
\end{prop}

\begin{proof}
To determine the boundary of $\Omega^{}_{AB}$, we choose special points
$P^{}_{i}$ in $\Omega^{}_{AB}$, see Figure~\ref{fig:win1} (for illustration)
and Table~\ref{tab:points} (for details)\footnote
{
We use the two dimensional geometry of the internal space here
explicitly, and, instead of going into cumbersome notations and explanations,
show some figures to clarify and assist the proofs.
}. 
We first show how these points are
determined. Demanding 
\begin{equation*}
P^{}_{2} = f^{}_{3}(P^{}_{1}), \qquad P^{}_{3} = f^{}_{1}(P^{}_{2}), \qquad 
P^{}_{4} = f^{}_{1}(P^{}_{3}), \quad \text{and} \quad P^{}_{1} =
f^{}_{3}(P^{}_{4})  
\end{equation*}
one gets the following fixed point equation
\begin{equation}\label{eq:P1}
P^{}_{1} = f^{}_{3} \circ f^{}_{1} \circ f^{}_{1} \circ f^{}_{3}(P^{}_{1}) =
\beta^4 \, P^{}_{1} + 6 \, \beta^2 + 2  
\end{equation}
for $P^{}_{1}$, and similar results hold for $P^{}_{2}$, $P^{}_{3}$ and
$P^{}_{4}$. The unique solution of~(\ref{eq:P1}) is
\begin{equation*}
P^{}_{1} = \frac12 (\beta^2-3\,\beta-1).
\end{equation*}
Choosing
\begin{equation*}
\begin{array}{lll}
P^{}_{5} = \frac12 (P^{}_{1}+P^{}_{3}), & P^{}_{6} = f^{}_{3}(P^{}_{3}),
& P^{}_{7} = f^{}_{1}(P^{}_{4}), \\ P^{}_{8} = f^{}_{2}(P^{}_{3}), &
P^{}_{9} = f^{}_{1}(P^{}_{1}),  & P^{}_{10} =
\frac12 (P^{}_{2}+P^{}_{3})
\end{array}
\end{equation*}
and setting $\tau$ to be the inversion in the center $P^{}_5$ ($\tau:\, z
\mapsto -z-\beta$) and $\kappa$ the one in the center $P^{}_{10}$
($\kappa:\, z \mapsto -z-\beta+1$), one can verify the following equations: 
\begin{equation*}
\begin{array}{lclclclcl}
P^{}_{1} &=& \tau(P^{}_{3}) \\
P^{}_{2} &=& \tau(P^{}_{4}) &=& \kappa(P^{}_{3}) \\
P^{}_{5} &=& \frac12 (P^{}_{2}+P^{}_{4}) &=& \frac12 (P^{}_{6}+P^{}_{9})
&=& \frac12 (P^{}_{7}+P^{}_{8}) \\
P^{}_{6} &=& \tau(P^{}_{9}) \\
P^{}_{7} &=& f^{}_{2}(P^{}_{1}) &=& f^{}_{3}(P^{}_{9}) &=& \tau(P^{}_{8}) \\
P^{}_{8} &=& f^{}_{1}(P^{}_{6}) &=& f^{}_{3}(P^{}_{2}) &=& \tau(P^{}_{8}) &=&
\kappa(P^{}_{9}) \\
P^{}_{10} &=& \frac12 (P^{}_{8}+P^{}_{9}).
\end{array}
\end{equation*}
For the mappings, one finds
\begin{equation*}
f^{}_{1} \circ \tau = \tau \circ f^{}_{3}, \qquad f^{}_{2} \circ \tau = \tau
\circ f^{}_{2} \quad \text{and} \quad f^{}_{3} \circ \tau = \tau \circ
f^{}_{1},
\end{equation*}
showing that $\Omega^{}_{AB}$ is inversion symmetric in the center
$P^{}_{5}$, i.e., $\tau(\Omega^{}_{AB}) =\Omega^{}_{AB}$.  

Denoting by $[P^{}_{2}, P^{}_{3}]$ the ``boundary'' between $P^{}_{2}$ and
$P^{}_{3}$ (the ``right edge''), one finds
\begin{equation*}
\begin{array}{lcl}
[P^{}_{2}, P^{}_{8}] &=& f^{}_{3} \circ \tau \circ f^{}_{1}([P^{}_{2},
P^{}_{3}]) \\[0pt]
[P^{}_{8}, P^{}_{9}] &=& f^{}_{1} \circ \tau \circ f^{}_{1} \circ \tau \circ
f^{}_{1}([P^{}_{2}, P^{}_{3}]) \\[0pt]
[P^{}_{9}, P^{}_{3}] &=& f^{}_{1} \circ \tau \circ f^{}_{1}([P^{}_{2},
P^{}_{3}]) 
\end{array}
\end{equation*}
and therefore the following IFS for $[P^{}_{2}, P^{}_{3}]$:
\begin{equation}\label{eq:IFS_boundary}
\begin{array}{lclclclcl}
[P^{}_{2}, P^{}_{3}] &=& g^{}_{1}([P^{}_{2}, P^{}_{3}]) &\cup&
g^{}_{2}([P^{}_{2}, P^{}_{3}]) &\cup& g^{}_{3}([P^{}_{2}, P^{}_{3}]), 
\end{array}
\end{equation}
where
\begin{equation}\label{eq:IFS_gs}
g^{}_1(z) = -\beta^2 \, z - \beta, \qquad g^{}_2(z) = (2 \beta^2+1) \, z +
\beta^2+1, \quad \text{and} \quad g^{}_3(z)=-\beta^2 \, z-\beta^2.
\end{equation}
Of course, we have not shown yet that $[P^{}_{2}, P^{}_{3}]$ 
really is (a piece of) the boundary of $\Omega^{}_{AB}$, so we just define
$[P^{}_{2}, P^{}_{3}]$ to be the unique compact solution of the IFS
(\ref{eq:IFS_boundary}), which is inversion symmetric in the center
$P^{}_{10}$ because  
\begin{equation*}
g^{}_{1} \circ \kappa = \kappa \circ g^{}_{3}, \qquad g^{}_{2} \circ \kappa =
\kappa \circ g^{}_{2} \quad \text{and} \quad g^{}_{3} \circ \kappa = \kappa
\circ g^{}_{1}.
\end{equation*}
Also, we know that $[P^{}_{2}, P^{}_{3}]$ is connected since we can start the
iteration with the straight line from $P^{}_{2}$ to $P^{}_{3}$. In each
iteration, the image remains a (piecewise smooth) path from $P^{}_{2}$ to
$P^{}_{3}$. 

With the mappings $f^{}_1$ and $\tau$, we get a boundary 
\begin{equation*}
\begin{array}{lclclclc}
[P^{}_{2}, P^{}_{3}] &\cup& f^{}_1([P^{}_{2}, P^{}_{3}]) &\cup& \tau([P^{}_{2},
P^{}_{3}] &\cup& \tau \circ f^{}_1([P^{}_{2}, P^{}_{3}]) &=\\ \hspace{0ex}
[P^{}_{2}, P^{}_{3}] &\cup& [P^{}_{3}, P^{}_{4}] &\cup& [P^{}_{4}, P^{}_{1}]
&\cup& [P^{}_{1}, P^{}_{2}]
\end{array}
\end{equation*}
around a simply connected open set $U$ (we will prove in the next proposition
that this boundary is non-self-intersecting).
Now one can show that only the boundaries of $\overline{f^{}_{i}(U)}$
intersect. Consider, for example, the region between
$f^{}_{2}(U)$ and $f^{}_{3}(U)$. Then the boundary
$[P^{}_{7}, P^{}_{8}]$ on $f^{}_{2}(U)$ is given by
\begin{equation*}
\begin{array}{lclcl}
[P^{}_{7}, P^{}_{8}] &=& f^{}_{2} \circ \tau \circ f^{}_{1}([P^{}_{2},
P^{}_{3}]) &\cup& f^{}_{2}([P^{}_{2}, P^{}_{3}]), 
\end{array}
\end{equation*}
while the one on $f^{}_{3}(U)$ is given by (taking orientation
into account)
\begin{equation*}
\begin{array}{lclcl}
[P^{}_{7}, P^{}_{8}] &=& f^{}_{3} \circ \kappa \circ g^{}_{2}([P^{}_{2},
P^{}_{3}]) &\cup& f^{}_{3} \circ g^{}_{3}([P^{}_{2}, P^{}_{3}]). 
\end{array}
\end{equation*}
It is easy to verify that
\begin{equation*}
f^{}_{2} \circ \tau \circ f^{}_{1} = f^{}_{3} \circ \kappa \circ g^{}_{2}
\quad \text{and} \quad f^{}_{2} = f^{}_{3} \circ \kappa \circ g^{}_{3}.
\end{equation*}
So the boundaries coincide. Similarly, one can check the region between
$f^{}_{1}(U)$ and $f^{}_{2}(U)$ (note that
$P^{}_{7}$ and $P^{}_{8}$ belong to all three sets
$\overline{f^{}_{i}(U)}$), and that the boundary of $U$ coincides with pieces
of the boundaries of the $f^{}_{i}(U)$ -- so the situation is as expected from
Figure~\ref{fig:win1}. Therefore, the IFS is just touching. Also, we now know
that $[P^{}_2,P^{}_3]$ is really a piece of the boundary of
$\Omega^{}_{AB}$.   
\end{proof}

\begin{prop}\label{prop:nonself_inter}
Let $\Omega^{}_{AB}$ be the unique compact solution of the
IFS $(\ref{eq:IFS-AB})$. Then its boundary is non-self-intersecting.
\end{prop}

\begin{figure}[t]
\begin{minipage}[t]{0.31\textwidth}
\centerline{\epsfxsize=\textwidth\epsfbox{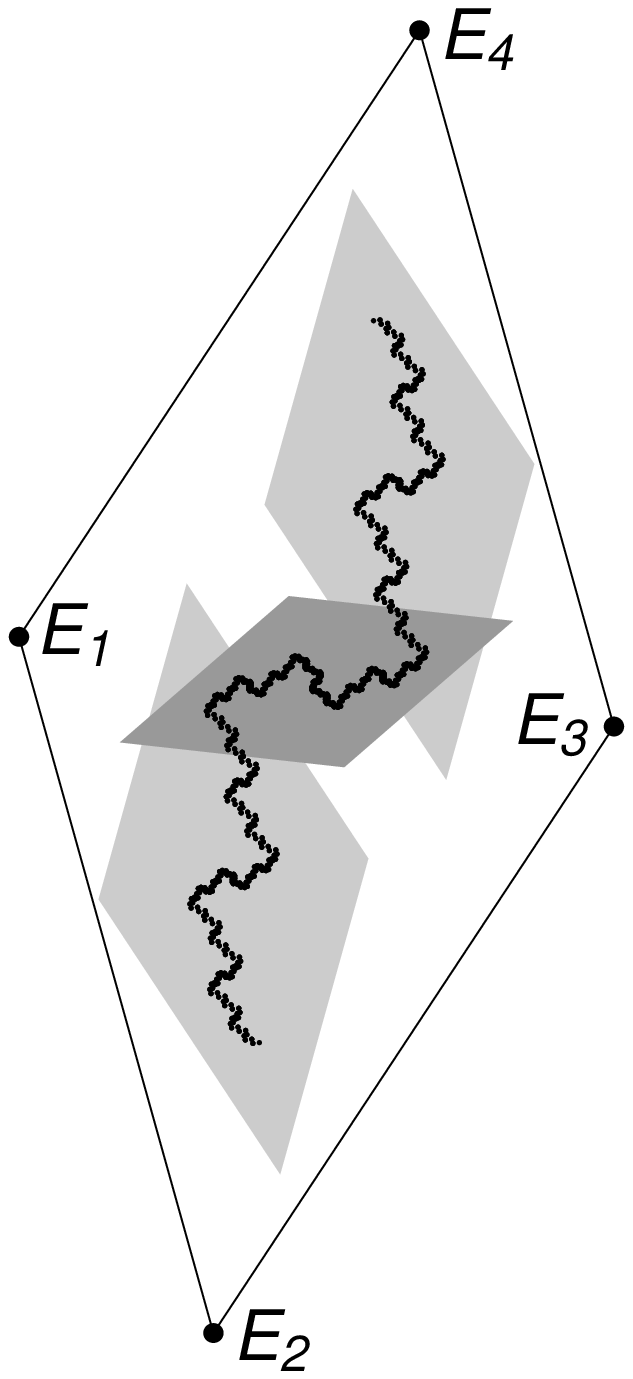}}
\setcaptionwidth{\textwidth}
\caption{\sloppy A rhombus $R$ with corners $E^{}_1,\ldots,E^{}_4$, such that
  $g^{}_{i}(R)\subset R$ for all $i$ (see~(\ref{eq:nsi_rhombus})). The
  $g^{}_{i}(R)$ are shaded in gray. For the
  coordinates of $E^{}_1,\ldots,E^{}_4$, see Table~\ref{tab:nsi_coord}.
  \label{fig:nsi_1}} 
\end{minipage}\hfill
\begin{minipage}[t]{0.31\textwidth}
\centerline{\epsfxsize=\textwidth\epsfbox{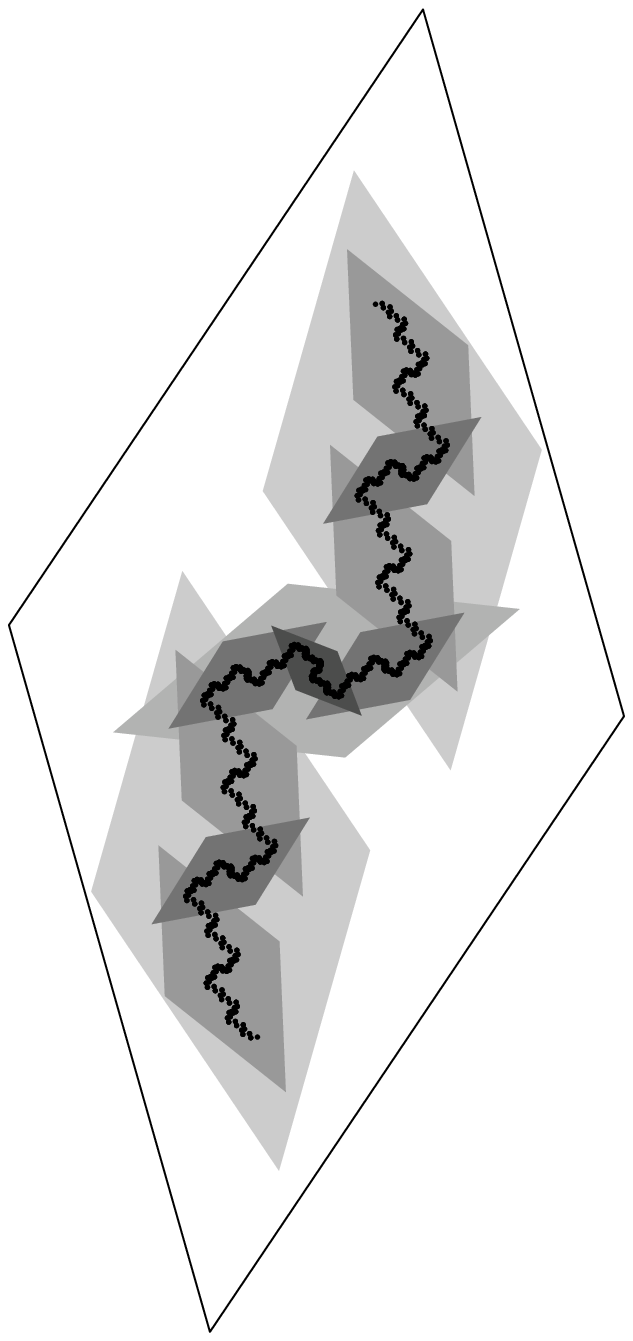}}
\setcaptionwidth{\textwidth}
\caption{\sloppy In addition to Figure~\ref{fig:nsi_1}, the second iteration
  of the rhombus $R$ is also shown, shaded in the dark
  gray. (\ref{eq:nsi_2stage}) is satisfied.\label{fig:nsi_2}} 
\end{minipage}\hfill
\begin{minipage}[t]{0.31\textwidth}
\centerline{\epsfxsize=\textwidth\epsfbox{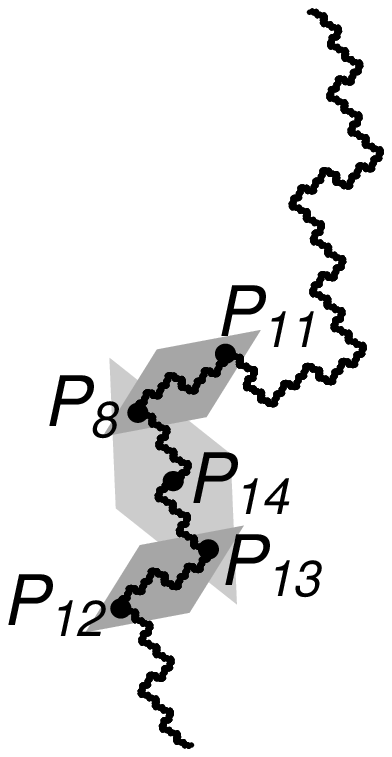}}
\setcaptionwidth{\textwidth}
\caption{\sloppy The situation around the ``joint'' $P^{}_{8}$. We have an
  inversion symmetry in the center $P^{}_{14}$. For the coordinates, see
  Table~\ref{tab:nsi_coord}. \label{fig:nsi_3}} 
\end{minipage}
\end{figure}

\begin{figure}[t]
\begin{minipage}[c]{0.3\textwidth}
\centerline{\begin{tabular}{|c|c|}\hline
Point & Coordinate \\ \hline \hline \rule[-1.5ex]{0ex}{4ex}
$E^{}_1$ & $P^{}_{10}-i\, \frac25\, (\beta^2-2\, \beta)$ \\
\rule[-1.5ex]{0ex}{4ex} 
$E^{}_2$ & $P^{}_{2}+\frac25\, (\beta^2-2\, \beta)$ \\ \rule[-1.5ex]{0ex}{4ex}
$E^{}_3$ & $P^{}_{10}+i\, \frac25\, (\beta^2-2\, \beta)$ \\
\rule[-1.5ex]{0ex}{4ex} 
$E^{}_4$ & $P^{}_{3}-\frac25\, (\beta^2-2\, \beta)$ \\ \rule[-1.5ex]{0ex}{4ex}
$P^{}_8$ & $\frac12 \, (\beta^2 - \beta + 1)$ \\ \rule[-1.5ex]{0ex}{4ex}
$P^{}_{11}$ & $\frac12 \, (9\, \beta^2+\beta+5)$ \\ \rule[-1.5ex]{0ex}{4ex}
$P^{}_{12}$ & $\frac12\, (-3\, \beta^2-3\, \beta-1)$ \\
\rule[-1.5ex]{0ex}{4ex}  
$P^{}_{13}$ & $\frac12\, (5\, \beta^2-\beta+3)$ \\ \rule[-1.5ex]{0ex}{4ex}
$P^{}_{14}$ & $\frac12\, (3\, \beta^2-\beta+2)$ \\ \hline 
\end{tabular}}
\vspace{1.5ex}
\tabcaption{\sloppy Coordinates of the points used in Figures~\ref{fig:nsi_1}
  and~\ref{fig:nsi_3}.\label{tab:nsi_coord}}
\end{minipage}\hfill
\begin{minipage}[c]{0.5\textwidth}
\centerline{\epsfxsize=\textwidth\epsfbox{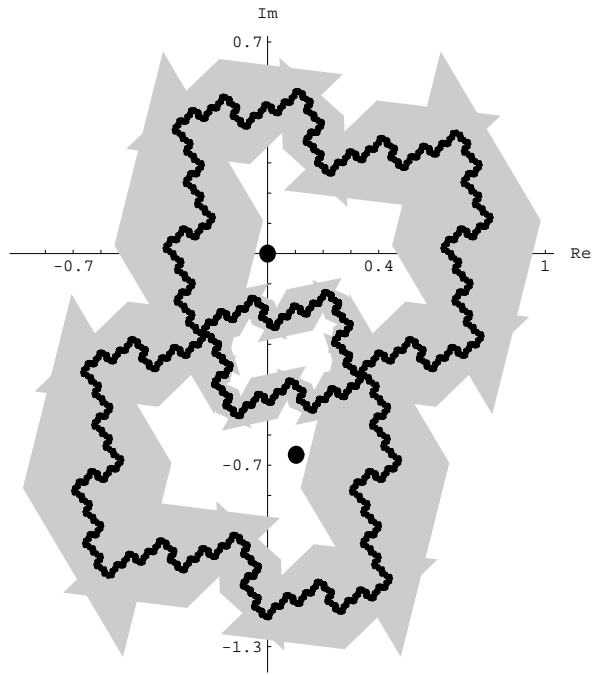}}
\setcaptionwidth{\textwidth}
\caption{\sloppy The points $0$ and $-\beta$ are represented as black
  dots. Also, the boundary is shown and the area of the rhombi (first
  iteration for all parts of the boundary) is shaded in gray. The boundary
  runs inside this shaded region, so $0$ and $-\beta$ are inner
  points. \label{fig:innerpoints}}  
\end{minipage}
\end{figure}

\begin{proof}
We specify a (closed) rhombus $R$ (which surrounds the boundary
$[P^{}_{2},P^{}_{3}]$ and the straight line from $P^{}_{2}$ to $P^{}_{3}$)
such that its iteration in~(\ref{eq:IFS-AB}) will not leave $R$, i.e., such
that 
\begin{equation}\label{eq:nsi_rhombus}
\left( g^{}_{1}(R) \cup g^{}_{2}(R) \cup g^{}_{3}(R) \right) \cap (\RR^{2}_{}
\setminus R ) = \varnothing.
\end{equation} 
Furthermore, we also require that
\begin{equation}\label{eq:nsi_g1g3}
g^{}_{1}(R) \cap g^{}_{3}(R) = \varnothing
\end{equation}
and
\begin{equation}\label{eq:nsi_2stage}
\begin{array}{ll}
g^{}_{1}(R) \cap g^{}_{2}(g^{}_{i}(R)) = \varnothing & \mbox{for } i\neq 1 \\
g^{}_{2}(R) \cap g^{}_{1}(g^{}_{i}(R)) = \varnothing & \mbox{for } i\neq 3 \\
g^{}_{2}(R) \cap g^{}_{3}(g^{}_{i}(R)) = \varnothing & \mbox{for } i\neq 1 \\
g^{}_{3}(R) \cap g^{}_{2}(g^{}_{i}(R)) = \varnothing & \mbox{for } i\neq 3.
\end{array}
\end{equation}
Such a rhombus $R$ exists, see~Figure~\ref{fig:nsi_1} for a picture of such a
rhombus that satisfies the conditions of~(\ref{eq:nsi_rhombus})
and~(\ref{eq:nsi_g1g3}) and Table~\ref{tab:nsi_coord} for the coordinates of
its corners (of course, we could also use a shape different from a
rhombus). This rhombus also satisfies~(\ref{eq:nsi_2stage}),
see~Figure~\ref{fig:nsi_2}.

Here, (\ref{eq:nsi_g1g3}) tells us that $[P^{}_{2},P^{}_{8}]$ and
$[P^{}_{9},P^{}_{3}]$ do not have a point in common; similar statements apply
for~(\ref{eq:nsi_2stage}). Each iterate of the rhombus is associated to a
corresponding iterate of the boundary $[P^{}_{2},P^{}_{3}]$ or the straight
line from $P^{}_{2}$ to $P^{}_{3}$ (denoted by $(P^{}_2,P^{}_3)$). We call two
rhombi at the same iteration level neighbouring if their corresponding
iteration of $(P^{}_{2},P^{}_{3})$ have a common endpoint. We see in 
Figure~\ref{fig:nsi_2} that only neighbouring rhombi intersect at the
second iteration level.  We show that for any iteration level only neighbouring
rhombi intersect.

We have verified the assertion for the first and second iteration level and
proceed inductively. Since $g^{}_{1}$, $g^{}_{2}$ and $g{}_{3}$ are affine, we
get the third iteration level as follows: The associated rhombi between
$P^{}_{2}$ and $P^{}_{8}$ are a scaled down (by $g^{}_{1}$) version of those
of the second level, therefore the assertion holds for them. Similarly for the
rhombi between $P^{}_{8}$ and $P^{}_{9}$ (by $g^{}_{2}$) and between
$P^{}_{9}$ and $P^{}_{3}$ (by $g^{}_{3}$). So the only critical points
remaining are the ``joints'' at $P^{}_{8}$ and $P^{}_{9}$. We show that at
these points also only neighbouring rhombi intersect and for this, we make use
of the self-similar structure of the boundary, see~Figure~\ref{fig:nsi_3}: The
boundary $[P^{}_{12},P^{}_{11}]$ is inversion symmetric in the center
$P^{}_{14}$. This is clear for $[P^{}_{13},P^{}_{8}] =
g^{}_{1}(g^{}_{3}([P^{}_{2},P^{}_{3}]))$ (and therefore $P^{}_{14} =
g^{}_{1}(g^{}_{3}(P^{}_{10}))$). But it also holds for $[P^{}_{12},P^{}_{13}]
= g^{}_{1}(g^{}_{2}([P^{}_{2},P^{}_{3}]))$ and $[P^{}_{8},P^{}_{11}] =
g^{}_{2}(g^{}_{1}(\kappa([P^{}_{2},P^{}_{3}])))$. So, since the assertion
holds around $P^{}_{13}$, it also holds around $P^{}_{8}$ by symmetry. Similar
arguments apply around $P^{}_{9}$. So the assertion holds for the third
iteration level, i.e., for the third iteration level only neighbouring rhombi
intersect. But the same argument applies to all further iteration levels. So
the assertion is true, i.e., for a given iteration level only neighbouring
rhombi intersect. Also note that each rhombus has two neighbouring rhombi
(with the exception of the ``starting'' and ``ending'' rhombi at $P^{}_{2}$
and $P^{}_{3}$ which only have one) and that there is no ``rhombus loop'',
i.e., going from $P^{}_{2}$ to $P^{}_{3}$ we cross each rhombus only once.

Now, suppose $[P^{}_{2},P^{}_{3}]$ is self-intersecting. Then there exist 
points $x,y \in [P^{}_{2},P^{}_{3}]$ ($x\neq y$) such that they are connected
in $[P^{}_{2},P^{}_{3}]$ in two different ways, $W^{}_{1}$ and $W^{}_{2}$ (and
we have a loop). We can choose points $u \in W^{}_{1}$ and  $v \in
W^{}_{2}$ such that $d(u,W^{}_{2}) = \min\{ d(v,z)\mathbin| z \in W^{}_{2} \}
>0$ ($W^{}_{2} \subset [P^{}_{2},P^{}_{3}]$ is compact) and
$d(v,W^{}_{1})>0$. But then, $u, v$ are in non-neighbouring rhombi for some
iteration level $N$ (and then for all iteration levels $n \ge N$), since the
length of a rhombus of the $N$th iteration level is at most $|\beta|_{}^{2N}
\cdot |E^{}_{2}-E^{}_{4}|$. So, we get a ``rhombus loop'' for this
iteration level by the rhombi which overlay $W^{}_{1}$ and $W^{}_{2}$. This is
a contradiction, therefore $[P^{}_{2},P^{}_{3}]$ is non-self-intersecting.    

From this single edge we proceed to all of the boundary. Here, 
critical are the ``joints'' $P^{}_{2}$, $P^{}_{3}$, etc., again, because
we get the other three parts by an affine map of this edge (e.g.,
$[P^{}_{3},P^{}_{4}]=f^{}_{1}([P^{}_{2},P^{}_{3}])$) and opposite edges (i.e.,
$[P^{}_{2},P^{}_{3}]$ and $[P^{}_{4},P^{}_{1}]$) do not overlap,
cf.~Figure~\ref{fig:innerpoints}. But at $P^{}_{2}$, an argument like the one
at $P^{}_{8}$ above applies, i.e., we have an inversion  symmetry of part of
the boundary in the center $\frac12\, (P^{}_{12}+P^{}_{13})$ (and similar for
the other ``joints''). This extends our findings to the entire boundary.
\end{proof}

This also implies that the boundaries of $\Omega^{}_{A}$, $\Omega^{}_{B}$ and
$\Omega^{}_{C}$, respectively their union $\Omega$, are
non-self-intersecting. Also, from the proof of the last proposition, we can
deduce the following.

\begin{coro}\label{coro:innerpoints}
The point $0$ is an inner point of $f^{}_{1}(\Omega^{}_{AB})\subset
\Omega^{}_{A}$ and $-\beta$ is an inner point of $f^{}_{3}(\Omega^{}_{AB})
\subset \Omega^{}_{B}$. 
\end{coro}

\begin{proof}
We again use the iteration of rhombi as in
Proposition~\ref{prop:nonself_inter} to show that the two points are really
inner points in the respective areas. For this, see
Figure~\ref{fig:innerpoints}, where the first iteration of the rhombi is
used for all parts of the boundary. Clearly, the points $0, -\beta$ are inner
points, which can easily be checked by a simple (though somewhat tedious)
calculation of distances.
\end{proof}

\begin{prop}\label{prop:setAB} 
Let $\Omega^{}_{AB}$ be the unique compact solution of the
IFS $(\ref{eq:IFS-AB})$. 
\begin{enumerate}
\item\label{en:statem0} $\Omega^{}_{AB}$ is inversion symmetric in the
  center $P^{}_5 = - \frac12 \beta$.
\item\label{en:statem1} $\Omega^{}_{AB}$ has Hausdorff dimension
  $h(\Omega^{}_{AB})=2$. 
\item\label{en:statem3} $\Omega^{}_{AB}$ has positive (Hausdorff and
  Lebesgue) measure (area).
\item\label{en:statem2} The boundary $\partial\Omega^{}_{AB}$ has vanishing
  (Lebesgue) measure.
\item\label{en:statem4} There is a periodic tiling of the plane with 
  $\Omega^{}_{AB}$ as prototile.
\end{enumerate}
\end{prop}

\begin{proof}
\ref{en:statem0} See proof of Proposition~\ref{prop:jt}.\\
\ref{en:statem1} Just touching implies the OSC, therefore
$h(\Omega^{}_{AB})=s(\Omega^{}_{AB})=2$. \\
\ref{en:statem3} The OSC for $\Omega^{}_{AB}$ (or any self-similar set with
similarity dimension $s$) is equivalent to the positive Hausdorff measure
condition $\mu^{s}_{}(\Omega^{}_{AB})>0$, where $\mu^{s}_{}$ denotes the
$s$-dimensional Hausdorff measure, see~\cite{Ban97} and references
therein. For Euclidean dimensions, Hausdorff and Lebesgue measure are
connected by a nonzero multiplicative constant.\\   
\ref{en:statem2} The similarity dimension
$\tilde{s}=s(\partial\Omega^{}_{AB})$ of the boundary is the solution of
(contraction constants given in (\ref{eq:IFS_gs})) 
\begin{equation*}
2 \, |\beta|^{2 \tilde{s}} + |\beta|^{3 \tilde{s}} = 1,
\end{equation*}
which is $\tilde{s} = -\log(\tau)/\log(|\beta|) \approx 1.22$ (where $\tau =
\frac12(1+\sqrt{5})$ is the golden ratio; the previous equation is solved by
$|\beta|^{-\tilde{s}} = \tau$). Therefore, the statement follows
from $h(\partial\Omega^{}_{AB}) \le s(\partial\Omega^{}_{AB})$.\\
\ref{en:statem4} Because of the inversion symmetries $\tau$ of
$\Omega^{}_{AB}$ and $\kappa$ of $\partial\Omega^{}_{AB}$ from $P^{}_2$
to $P^{}_3$ (see proof of Proposition~\ref{prop:jt}), the ``right edge'' 
and the ``left edge'' differ only by a translation $P^{}_{2} - P^{}_{1}=1$,
and, similarly, the ``upper edge'' and the ``lower edge'' differ by $P^{}_{2}
-  P^{}_{3}=\beta^2-2\,\beta$.  
\end{proof} 

\begin{prop}\label{prop:lattice}
Let $\Omega^{}_{A},\, \Omega^{}_{B},\, \Omega^{}_{C}$ be the solution
of the IFS $(\ref{eq:IFS})$, and $\Omega = \Omega^{}_{A}\cup
\Omega^{}_{B}\cup \Omega^{}_{C}$. Then $\Omega$ is a compact set,
homeomorphic to a disc, with positive area. The boundary $\partial\Omega$
is a fractal of vanishing Lebesgue measure, which is
non-self-intersecting. The set $\Omega$ admits a lattice tiling of $\RR^2$,
where the lattice is spanned by $P^{}_{2} - P^{}_{6}=-\beta+1$ and $P^{}_{2} -
P^{}_{3}= \beta^2-2\,\beta$.  
\end{prop}

\begin{proof}
It is clear from our construction that $\Omega$ is a compact set with simply
connected interior. We have also seen that the boundary is connected and
consists of finitely many pieces, each of which is obtained from a
construction as used in the proof of Proposition~\ref{prop:jt}. So, $\Omega$
must be homeomorphic to a disc.
The remaining statements follow directly from
Propositions~\ref{prop:nonself_inter} and~\ref{prop:setAB}, because the
mappings in (\ref{eq:affine}) are 
affine and the just touching property also holds for
$\Omega=\Omega^{}_{A}\cup\Omega^{}_{B}\cup\Omega^{}_{C}$. Since we also know
the boundary of $\Omega$ (we have an IFS for every part of it), we can also
verify the translation vectors by comparing the corresponding iterated
function systems. Also, see Figure~\ref{fig:win2} for a depiction of these
vectors.   
\end{proof}

\begin{coro}
$\varLambda(\Omega)$ is a regular model set.\qed
\end{coro}

We can calculate the volume\footnote
{
Note that the discriminant of $\QQ(\alpha)$ is $-59$. The volume
$|\varGamma|$ is proportional to the square root of the absolute value of the
discriminant. The proportional constant is one factor of $\frac12$ because
there is one complex conjugate pair $\beta, \overline{\beta}$ of algebraic
conjugates of $\alpha$, see~\cite[Chapter II, Section 4.2, Theorem
2]{BS66}. Note that we also have a formula for $|\IM(\beta)|$ in terms of
$\alpha$ by this: 
\begin{equation*}
|\IM(\beta)| = \frac1{2\cdot\sqrt{59}}\,(-8\,\alpha^2+25\,\alpha-6).
\end{equation*}
} 
\begin{equation}\label{eq:vol}
|\varGamma| = |\det( \bs{v}^{}_{A}, \bs{v}^{}_{B}, \bs{v}^{}_{C})| =
 |\IM(\beta)|(3\,\alpha^2-4\,\alpha) = \frac12 \sqrt{59}\approx 3.84 
\end{equation}
of the fundamental domain of $\varGamma$.
And because of the periodic tiling of the plane with $\Omega$ as a prototile,
it is also easy to calculate the area $\mu^{}_{\textnormal{int}}(\Omega)$ of
$\Omega$: as $\mu^{}_{\textnormal{int}}(\partial\Omega) = 0$,
$\mu^{}_{\textnormal{int}}(\Omega)$ equals the area of a fundamental domain
of the corresponding lattice of periods. This gives 
\begin{equation}\label{eq:area}
\mu^{}_{\textnormal{int}}(\Omega) = |\det((P^{}_{2} - P^{}_{6}),(P^{}_{2} -
P^{}_{3}))| = |\IM(\beta)|(\alpha^2-\alpha) = \frac12
\frac{1}{\sqrt{59}}(3\,\alpha^2-2\,\alpha+17)\approx 1.77. 
\end{equation}
Then the following lemma applies.

\begin{lem} 
Let $\varGamma$ be a lattice in $\RR \times \RR^m$, $|\varGamma|$ be the
volume of a measurable fundamental domain of $\varGamma$ in $\RR \times
\RR^m$ with respect to the product $\mu \otimes \mu^{}_{\textnormal{int}}$ of
the Lebesgue measures $\mu$, $\mu^{}_{\textnormal{int}}$ on $\RR$, $\RR^m$,
respectively. Assume that we have a cut-and-project scheme like in
$(\ref{eq:cutproj})$. If $\Omega$ is a bounded subset of $\RR^m$ with almost
no boundary, then the density $\operatorname{dens}(\varLambda(\Omega))$ of
the corresponding regular model set in $\RR$ is 
\begin{equation*}
\operatorname{dens}(\varLambda(\Omega))=
\frac{\mu^{}_{\textnormal{int}}(\Omega)}{|\varGamma|}. 
\end{equation*}
\end{lem} 

\begin{proof}
This follows from \cite[Proposition 2.1]{Sch98} because the projection $\pi$
is one-to-one on $\varGamma$ by construction.
\end{proof}

With (\ref{eq:avlen}), (\ref{eq:vol}) and (\ref{eq:area}), it is now easy to
check that the density of the model set $\varLambda(\Omega)$ and the density
of $\CPS$ are equal, i.e.,
\begin{equation}\label{eq:dens}
\frac{\mu^{}_{\textnormal{int}}(\Omega)}{|\varGamma|} = \frac1{\ell}.
\end{equation} 

\begin{prop}\label{prop:density}
The sequence $\CPS$ is a subset of $\varLambda(\Omega)$. Further,
they differ at most on positions of zero density and therefore have the same
pure point diffraction spectrum.
\end{prop}

\begin{proof}
We choose $\bs{0},-\bs{v}^{}_{B} \in \varGamma$. Then their projections into
internal space are elements of the attractor $\Omega$, because
$f^{}_1(0) = 0$ and $f^{}_3(-\beta) = -\beta$. But starting with these two
points, the iteration of the IFS in internal space just corresponds to the
iteration in (\ref{eq:BA}), respectively ~(\ref{eq:IFS_physical}), in physical
space. Therefore, $\CPS \, \subset \, \varLambda(\Omega)$, because the star
map of all iterates of $0$ and $-\alpha$ (i.e., $0$ and $-\beta$) stay in
$\Omega$. Equation~(\ref{eq:dens}) shows that both sequences have the same
density. So they can at most differ on positions of zero density.

Regular model sets have a pure point diffraction
spectrum, see~\cite{BD00,BM02,Sch00} and references therein. Therefore, the
diffraction spectrum of $\varLambda(\Omega)$ is pure point, and $\CPS$,
differing at most at positions of zero density, has the same spectrum by an
argument in~\cite{Hof95}. 
\end{proof}

\begin{theorem}\label{thm:regular}
$\CPS$ is a regular model set (except possibly for positions of
zero density) and has a pure point diffraction
spectrum. Its autocorrelation is a norm almost periodic point measure,
supported on a uniformly discrete subset of $\ZZ[\alpha]$.
\end{theorem}

\begin{proof}
The first assertion follows from Proposition~\ref{prop:density}. The
autocorrelation measure (see \cite{Hof95} for details) is supported on
$\varLambda(\Omega) - \varLambda(\Omega) \subset \ZZ[\alpha]$. Since
$\varLambda(\Omega)$ is a model set, it is also a Meyer set, hence
$\varLambda(\Omega) - \varLambda(\Omega)$ is Delone. The norm almost
periodicity follows from \cite[Theorem 5]{BM02}. 
\end{proof}

What we have proved so far is enough to calculate the diffraction spectrum of
$\CPS$ and Kol$(3,1)$, see Section~\ref{sec:diff}. But in the next section, we
want to show that $\varLambda(\Omega)$ really \textit{equals} $\CPS$. 

\smallskip
\noindent
\textsc{Remarks}: For unimodular substitutions of Pisot type, i.e.,
Kol$(2m\pm 1,2m\mp 1)$ with $m \ge 1$, the procedure is essentially the
same. Unfortunately, the IFS does not decouple like in (\ref{eq:IFS-AB}) for
$m>1$, which makes it technically more involved. For non-unimodular
substitutions of Pisot type, the internal space is more complicated in having
additional $p$-adic type components, see~\cite{Sie02, GK97, BMS98, LM01, LMS}
for further details and examples.  

\smallskip
\noindent
The sequences $\Sigma\textnormal{Kol}(p,q)$, which are not of Pisot type,
do not have a pure point spectral component outside $k=0$ by an argument
in~\cite{BT86}, see~\cite{Diplom} for details. 

\smallskip
\noindent
Some of the results given have been studied extensively under the name of
``Rauzy fractal'', e.g., that the windows have non-empty interior, that the
windows do not overlap in this case (this follows from the so-called strong
coincidence condition) and also the periodic tilability seems to follow from
results in~\cite{AI01, CS01, SW03}. But we also need the lattice of the
periodic tiling explicitly, as well as the induced IFS~(\ref{eq:IFS_boundary})
for the boundary. Therefore, we opted to give an elementary and complete
derivation here.  

\section{$\CPS$ is a generic model set}

We set $\diamondset{\Omega} = \openset{\Omega}^{}_{A} \cup
\openset{\Omega}^{}_{B} \cup \openset{\Omega}^{}_{C} \subset \Omega$
and $\partial \diamondset{\Omega} = \Omega \setminus
\diamondset{\Omega} = \partial\Omega^{}_{A} \cup
\partial\Omega^{}_{B} \cup \partial\Omega^{}_{C}$. Then we can improve a
statement of Proposition~\ref{prop:density}.    

\begin{prop} \label{prop:inside}
$\CPS$ is equal to the model set
  $\varLambda(\diamondset{\Omega})$.  
\end{prop}

\begin{proof}
$\CPS \subset \varLambda(\diamondset{\Omega})$:
Note that the mappings $f^{}_i$ ($i \in \{0,1,3\}$) of (\ref{eq:IFS}') are
similarities (all directions are contracted by the same factor, here
$|\beta|$). Therefore, they map balls around $x$ to balls around 
$f^{}_i(x)$. Furthermore, they map balls in 
$\openset{\Omega}^{}_{k}$ to balls in $\openset{\Omega}^{}_{\ell}$
($k,\ell \in \{A,B,C\}$). Since the starting points of
the iteration~(\ref{eq:BA}) in internal space, namely $-\beta$ and $0$, are
inner points of $\Omega^{}_{B}$ and $\Omega^{}_{A}$ by
Corollary~\ref{coro:innerpoints}, one can also find 
balls of radius $\varepsilon > 0$ around $-\beta$ and $0$ which lie entirely in
$\openset{\Omega}^{}_{B}$ and $\openset{\Omega}^{}_{A}$,
respectively. Since the iteration in physical space corresponds to the IFS in
internal space, the star map of an arbitrary point in $\CPS$ is thus a point of
$\diamondset{\Omega}$. \\
$\varLambda(\diamondset{\Omega}) \subset \CPS$:
Suppose $x \in \varLambda(\diamondset{\Omega}) \setminus
\CPS$. Then $x^{\star}$ and all its iterates of the
IFS~(\ref{eq:IFS}) are in $\diamondset{\Omega}$, by the same reasoning as
before. Furthermore, the mappings $f^{}_i$ ($i \in \{0,1,3\}$) are affine
similarities and therefore all iterates of $x^{\star}$ are disjoint to all of
the iterates of $-\beta$ and $0$. But then
\begin{equation*}
\operatorname{dens} \CPS < \operatorname{dens}
\varLambda(\diamondset{\Omega}), 
\end{equation*}
because the set of iterates of $x$ under inflation has positive density in
$\varLambda(\diamondset{\Omega})$. This contradicts Proposition
\ref{prop:density}. 
\end{proof}

Note that not only the original sequence with $3$'s and $1$'s is inversion
symmetric (see~(\ref{eq:31})), but also the positions in
$\varLambda(\Omega^{}_{AB})$. 

\begin{coro}
$\varLambda(\Omega^{}_{AB})$ is inversion symmetric in the center $-\frac12
\alpha$. 
\end{coro}

\begin{proof}
The starting points $-\beta$ and $0$ and the IFS~(\ref{eq:IFS-AB}) are
inversion symmetric in the center $P^{}_5$. This corresponds, in physical space
(by Galois conjugation), to inversion symmetry in the center $- \frac12
\alpha$. 
\end{proof}

For the following, we need some more definitions, see~\cite{Baa02}
and~\cite{Moo00}. If $\varLambda$ is a discrete point set in $\RR^d$, we call 
$\mathcal{P}^{}_{r}(u)$ an \textit{$r$-patch} of a point $u \in \varLambda$, if
$\mathcal{P}^{}_{r}(u) = \varLambda \cap B^{}_{r}(u)$, where $B^{}_{r}(u)$ is
the ball of radius $r$ about $u$. Often, we are only interested in the set of
$\{\mathcal{P}^{}_{r}(u) \mathbin| r>0, u \in \varLambda \}$ and an element of
this set is simply called a \textit{patch}. Two structures
$\varLambda^{}_1$ and $\varLambda^{}_2$ are \textit{locally indistinguishable}
(or \textit{locally isomorphic} or \textit{LI}) if each patch of
$\varLambda^{}_1$ is, up to translation, also a patch of $\varLambda^{}_2$ and
vice versa. The corresponding equivalence class is called
\textit{LI-class}. 

A discrete structure $\varLambda$ is \textit{repetitive}, if for every $r>0$
there is a radius $R(r)>0$ such that within each ball of radius $R(r)$, no
matter its position in $\RR^d$, there is at least one translate of each
$r$-patch. Note that every primitive substitution generates a repetitive
sequence, wherefore $\CPS$ is repetitive.

We now look at generic model sets and show that $\varLambda(\Omega)$
is actually generic.

\begin{lem} \label{lem:C_dense}
The set $C = \{ c \in \RR^2 \mathbin| (c + \partial\diamondset{\Omega}) \,
\cap \, L^{\star} = \varnothing \}$ is dense in $\RR^2$, especially
$c+L^{\star} \subset C$ for $c \in C$.
\end{lem}

\begin{proof}
The set $C$ is not the empty set by Baire's category theorem
($\partial\diamondset{\Omega}$ is a meager set, $L^{\star}$ is countable), by
standard arguments, which in this context first appeared in~\cite[Section
2.2.2]{Sch93}, also see~\cite{BMS98}.
But if $c \in C$, then $c+t \in C,\; \forall t \in L^{\star}$ (notice that
$L^{\star}$ is an Abelian group), and $L^{\star}$ is dense. 
\end{proof}

\begin{prop}\label{prop:LI}\makebox[0ex]\\
\begin{enumerate}
\item\label{en:C_generic}
The model set $\varLambda(c+\Omega)$ is repetitive and generic for $c \in
C$. 
\item\label{en:C_LI}
The model sets $\varLambda(c+\Omega)$ and $\varLambda(\tilde{c}+\Omega)$
are LI for $c, \tilde{c} \in C$.
\item\label{en:LI} $\CPS$ and $\varLambda(c+\Omega)$ are LI for $c \in C$.
\end{enumerate}
\end{prop}

\begin{proof}
\ref{en:C_generic} The model set $\varLambda(c+\Omega)$ is generic by the
definition of the set $C$. It is repetitive by~\cite[Theorem 6]{Sch98}
and~\cite[Proposition 3.1]{Sch00}.\\  
\ref{en:C_LI} This is a by now standard argument, apparently first used
in~\cite[Lemma 2.1]{Sch93}.\\ 
\ref{en:LI} Since $\varLambda(c+\Omega)$ is repetitive, by~\cite[Lemma
1.2]{Sch93} it is enough to check that every patch of
$\CPS$ also occurs\footnote
{
That every patch of $\varLambda(c+\Omega)$ is also one of $\CPS$, follows then
together with the repetivity of $\varLambda(c+\Omega)$.
}
as a patch of $\varLambda(c+\Omega)$: Let
$\mathcal{P}$ be a patch of $\CPS$. Then $\mathcal{P}^{\star} \subset
\diamondset{\Omega}$ by Proposition~\ref{prop:inside}, and since $\mathcal{P}$
is a finite patch, we even know that there is an $\varepsilon > 0$ such that
$\mathcal{P}^{\star} \subset \diamondset{\Omega}^{}_{\varepsilon}$, where
$\diamondset{\Omega}^{}_{\varepsilon} = \{t \in \RR^2 \mathbin| t \in
\diamondset{\Omega}, \,
\operatorname{dist}(t,\partial\diamondset{\Omega}) > \varepsilon \}$. The
set $C$ is dense by Lemma~\ref{lem:C_dense}, therefore there is a $\tilde{c}
\in C$ such that $\mathcal{P}^{\star} \subset \tilde{c}+\Omega$. Then
$\mathcal{P} \subset \varLambda(\tilde{c}+\Omega)$ and $\CPS$ and
$\varLambda(\tilde{c}+\Omega)$ are LI. Since LI is an equivalence relation,
it follows from~\ref{en:C_LI} that $\CPS$ and $\varLambda(c+\Omega)$ are
LI for every $c \in C$. 
\end{proof}

\begin{prop}\label{prop:C0}
Define $C^{}_0 = \{ c \in C \mathbin| 0 \in c+\Omega \}$. Then, for every
$r>0$ and $c \in C^{}_0$, we get
\begin{equation*}
B^{}_{r}(0) \cap \CPS = B^{}_{r}(0) \cap
\left(\varLambda(c+\Omega)+t\right) =  B^{}_{r}(0) \cap
\varLambda(c+t^{\star}+\Omega) 
\end{equation*}
for an appropriate $t \in \ZZ[\alpha]$. Additionally, it follows that
$c+t^{\star} \in C^{}_0$. 
\end{prop} 

\begin{proof}
By Proposition~\ref{prop:LI}\ref{en:LI} we know that $\varLambda(c+\Omega)$
and $\CPS$ are LI, therefore the patch $B^{}_{r}(0) \cap \CPS$ occurs
somewhere in $\varLambda(c+\Omega)$. By the choice of $C^{}_0$, we can
translate this patch in $\varLambda(c+\Omega)$ with $t \in \ZZ[\alpha]$ to the
origin. In internal space, this is a translation $t^{\star}$ and by
Lemma~\ref{lem:C_dense} we have $c+t^{\star} \in C$. But $\CPS$ has a point at
the origin, so we even have $c+t^{\star} \in C^{}_0$.
\end{proof}

\begin{theorem}\label{thm:generic}
$\CPS$ is a regular generic model set.
\end{theorem}

\begin{proof}
Since $\CPS$ is the fixed point of a primitive substitution, it is repetitive,
and the corresponding dynamical system is minimal, see~\cite[Proposition
3.1]{Sch00}.  

By Proposition~\ref{prop:LI}\ref{en:LI} the model set
$\varLambda(c+\Omega)$ for $c \in C$ is in the LI-class of $\CPS$, hence the
latter must be the limit of some sequence $(t^{}_{i})^{}_{i}$ of translations
of $\varLambda(c+\Omega)$, where the translates $t^{}_{i}$ can be restricted
to elements of $\ZZ[\alpha]$ and $c+t^{}_{i} \in C^{}_{0}$ by
Proposition~\ref{prop:C0}.  

However, if $\partial\Omega\cap L^{\star}\neq\varnothing$, there is a point
$\tilde{x} \in \partial\diamondset{\Omega}\cap L^{\star}\cap\openset{\Omega}$
with $\tilde{x} \in \partial\Omega^{}_{A} \cap \partial\Omega^{}_{B}$ (so
$\tilde{x}$ lies on the common boundary of $\Omega^{}_{A}$ and
$\Omega^{}_{B}$). This is because by appropriate combinations of 
the mappings $g^{}_1, g^{}_2, g^{}_3, \tau, \kappa, f^{}_0, f^{}_1, f^{}_2,
f^{}_3, f^{}_4$ and the translations by $P^{}_2-P^{}_1, P^{}_2-P^{}_3,
P^{}_2-P^{}_6$ of Section~\ref{sec:IFS}, which all map $L^{\star}$ onto
$L^{\star}$, we can ``move'' points on $\partial\diamondset{\Omega}$ from
every ``edge'' to every other ``edge''\footnote
{
We even get that, with one point $x \in \partial\diamondset{\Omega} \cap
L^{\star}$, there is a dense set of points in $\partial\diamondset{\Omega}
\cap L^{\star}$, because we can always ``move'' $x$ to the edge
$[P^{}_2,P^{}_3]$, apply the IFS~(\ref{eq:IFS_boundary}) there and ``move''
this edge, with now dense points, to every other edge.
}. 
We have $\varepsilon^{}_0 = \operatorname{dist}(\tilde{x},
\partial\diamondset{\Omega} \setminus (\partial\Omega^{}_{A} \cap
\partial\Omega^{}_{B})) > 0$.  

The inverse star image of this point must then be in any limit of sequences
$\varLambda(t^{\star}_{i}+c+\Omega)$ with $c+t^{\star}_{i} \to 0$, but it
is not in $\CPS$ --- which is a contradiction. So no such point can exist and
$\partial\Omega\cap L^{\star}=\varnothing$. 

Regularity was established in Theorem~\ref{thm:regular} together with
Proposition~\ref{prop:inside}.  
\end{proof}

This argument is rather general and applies in other situations as well. For a
more elementary proof see the appendix.

\smallskip
\noindent
\textsc{Remark}: 
By Proposition~\ref{prop:lattice} we know that $\bigcup_{t \in G}
(t+\Omega) = \RR^2$, where $G=\langle-\beta+1,\, \beta^2-2 \, \beta
\rangle_{\ZZ}$  is a rank $2$ free Abelian group (a $2$-dimensional lattice)
and by Theorem~\ref{thm:generic} that $L^{\star} = \ZZ[\beta] \, \subset \,
\dot{\bigcup}_{t \in G} (t + \openset{\Omega})$,
where $\dot{\cup}$ denotes disjoint union. In physical space, we get
\begin{equation*}
L = \ZZ[\alpha] = \dot{\bigcup_{s \in G'}} (s + \CPS), 
\end{equation*}
where 
\begin{equation*}
G'=\langle-\alpha+1,\, \alpha^2-2 \, \alpha \rangle_{\ZZ} = \langle
\ell^{}_{C}-\ell^{}_{B}, \,\ell^{}_{A}-\ell^{}_{B}\rangle_{\ZZ},
\end{equation*}
i.e., $\ZZ[\alpha]$ is the disjoint union of translates of the regular,
generic model set $\CPS$. The set of translations needed is a rank $2$
subgroup, whose Galois dual in $\ZZ[\beta]$ is a lattice. But we can also write
\begin{equation*}
L=\ZZ[\alpha] = \dot{\bigcup_{r \in \CPS}} (r + G'),
\end{equation*} 
so $L/G'$ is a coset system with the structure of a model set. Now, let
$\lambda(m)$ be the $m$-th element of $\CPS$ ($m \in \ZZ$). Then one can show
that the induced group structure on this coset system is $(\lambda(m)) +
(\lambda(n)) = (\lambda(m+n))$, i.e., it is the action of $\ZZ$ on
$\CPS$. This group structure lines up with the deformation in the next
section.  

\section{Deformation and Diffraction \label{sec:diff}}

In the cut-and-project scheme~(\ref{eq:cutproj}) with $H = \RR^m$, let
$\varphi: \RR^{m} \to \RR^{d}$ be a continuous function with compact support
(e.g., $\Omega$). We call   
\begin{equation*}
\varLambda^{}_{\varphi}(\Omega) = \{ x + \varphi(x^{\star}) \mathbin|
x \in \varLambda(\Omega) \}
\end{equation*}
a \textit{deformed model set} if it is also a Delone set, see~\cite{BD00}. The
model set $\varLambda(\Omega)$ can be seen as 
deformed model set where the associated function $\varphi$ is trivial, i.e.,
$\varphi\equiv 0$. The diffraction spectrum of (deformed) model sets
(where each of its points is represented by a normalized Dirac measure, say)
can be calculated explicitly, see~\cite{BD00} for details. We write
$\delta^{}_{k}$ for the Dirac measure at $k$, i.e., $\delta^{}_{k}(f)=f(k)$
for $f$ continuous. Also, we need the dual of a lattice $\varGamma \subset
\RR^{n}$ defined as 
\begin{equation*}
\varGamma^{\ast} := \{ y \in \RR^{n} \mathbin| x.y \in \ZZ, \,\forall x \in
\varGamma \},
\end{equation*}
with $x.y$ denoting the Euclidean scalar product.
 
\begin{prop}\cite{BD00}
Let $\varLambda^{}_{\varphi}(\Omega)$ be a deformed model set in $\RR^d$
constructed with a regular model set $\varLambda(\Omega)$ and a continuous
function $\varphi$ of compact support. Then, the diffraction pattern of
$\varLambda^{}_{\varphi}(\Omega)$ is the positive pure point measure 
\begin{equation*}
\hat{\gamma} = \sum_{k \in
  \pi(\varGamma^{\ast})}|c^{}_{k}(\varLambda^{}_{\varphi}(\Omega))|^2\,
  \delta^{}_{k}, 
\end{equation*}
where $\varGamma^{\ast}$ is the dual lattice of $\varGamma$, $\delta^{}_{k}$
is the Dirac measure at $k$ and $c^{}_{k}(\varLambda^{}_{\varphi}(\Omega))$
is the Fourier-Bohr coefficient of $\varLambda^{}_{\varphi}(\Omega)$ at
$k$. This Fourier-Bohr coefficient exists and has the value  
\begin{equation} \label{eq:FBcoeff}
c^{}_{k}(\varLambda^{}_{\varphi}(\Omega)) = \left\{ \begin{array}{cl}
\frac1{|\varGamma|} \int_{\Omega} e^{-2\pi i \, (k.\varphi(y) -
  k^{\star}.y)} dy, & \text{if } (k,k^{\star}) \in \varGamma^{\ast}, \\  0, &
\text{otherwise.} \\
\end{array} \right. 
\end{equation} \qed
\end{prop}

\noindent
For a regular model set $\varLambda(\Omega)$ (where $\varphi\equiv 0$), the 
Fourier-Bohr coefficient is just given by the (inverse) Fourier transform of
the characteristic function of the window $\Omega$. 
   
For $\CPS$, we note that the support $F$ of the spectrum is dense
in $\RR$ since it is given by the $\ZZ$-span of the projection of the dual
lattice vectors, i.e.,
\begin{equation} \label{eq:supp}
F = \ZZ \cdot \pi(\bs{v}^{\ast}_{A}) + \ZZ\cdot \pi(\bs{v}^{\ast}_{B})
+ \ZZ\cdot \pi(\bs{v}^{\ast}_{C}),
\end{equation}
but $\pi(\bs{v}^{\ast}_{B}) = (\alpha-1)\, \pi(\bs{v}^{\ast}_{A})$ and
therefore they are linearly independent over $\QQ$.

To deform $\CPS$ to Kol$(3,1)$, we make the linear ansatz
$\varphi(x^{\star}) = a \,x^{\star}_{1} + b \, x^{\star}_{2}$, where
$x^{\star}_{i}$ denotes the $i$th Cartesian 
component of the vector $x^{\star} \in \RR^2$. With this $\varphi$, we now
deform all bond lengths $\ell^{}_{i}$ to the average bond length $\ell$, i.e.,
we have to solve the following linear system of equations ($i \in \{A,B,C\}$):
\begin{equation*}
(\bs{v}_{i})_1 + a \,(\bs{v}_{i})_2 + b \,(\bs{v}_{i})_3 = \ell.
\end{equation*}        
This over-determined system is solved by
\begin{equation*}
\begin{split}
& a = \ell - 1 = \frac12(-\alpha^2+\alpha+5)\approx 1.17 \quad \text{and} \\
& b = \frac{\IM(\beta)}{59}\cdot(-\alpha^2-17\,\alpha+31) =
\frac1{2\,\sqrt{59}^3} (-413\,\alpha^2+885\,\alpha-59) \approx - 0.13.
\end{split}
\end{equation*}
Due to the linearity of $\varphi$ (and the positivity of the bond lengths
involved), this deformation does not alter the order of
the points (i.e., for $x, x' \in \varLambda(\Omega)$ with
$x<x'$, we always have $x + \varphi(x^{\star})
<x' + \varphi({x'}^{\star})$). Note that the support $F$ of the spectrum
stays the same as in (\ref{eq:supp}), only the Fourier-Bohr coefficients
change.  

The positions in $\varLambda^{}_{\varphi}(\Omega^{}_{i})$ are now subsets of
$\ell \cdot \ZZ$. To be more precise, we even have
\begin{equation*}
\varLambda^{}_{\varphi}(\Omega^{}_{A}) \cup
\varLambda^{}_{\varphi}(\Omega^{}_{B}) \cup
\varLambda^{}_{\varphi}(\Omega^{}_{C}) \, = \, \ell \cdot \ZZ.
\end{equation*}    
Because of this embedding into $\ell \cdot \ZZ$, the diffraction spectrum of
each of the aperiodic sets $\varLambda^{}_{\varphi}(\Omega^{}_{i})$ is
$(\ell \cdot \ZZ)^{\ast}$-periodic~\cite{Baa01}, i.e., it is periodic with
period $1/\ell$ (note that $\ZZ^{\ast} = \ZZ$; the diffraction spectrum of
$\CPS$ is not periodic). This might not be obvious from
(\ref{eq:FBcoeff}) at first sight, but for $n \in \ZZ$ we have (note that
$\pi(\bs{v}^{\ast}_{i}) = \rho^{}_{i}/\ell$) 
\begin{equation}\label{eq:invlen}
\frac{n}{\ell} = n \,
\pi(\bs{v}^{\ast}_{A}+\bs{v}^{\ast}_{B}+\bs{v}^{\ast}_{C}) = \pi(n \,
\bs{v}^{\ast}_{A}+n \, \bs{v}^{\ast}_{B}+n\, \bs{v}^{\ast}_{C}), 
\end{equation}
therefore, for every $k'=k+\frac{n}{\ell}$ with $(k,k^{\star}) \in
\varGamma^{\ast}$, there is also a ${k'}^{\star}$ with $(k',{k'}^{\star}) \in
\varGamma^{\ast}$ given by
\begin{equation*}
{k'}^{\star}=k^{\star} + \pi_{\textnormal{int}}^{}(n \, \bs{v}^{\ast}_{A}+n \,
\bs{v}^{\ast}_{B}+n\, \bs{v}^{\ast}_{C}). 
\end{equation*}
But with the chosen $\varphi$ we get
\begin{multline}\label{eq:bra2ket}
(k'-k).\varphi(y) - ({k'}^{\star} - k^{\star}).y \\ = n \, y^{}_1
\left[\frac{a}{\ell} - \left( \pi_{\textnormal{int}}^{}(\bs{v}^{\ast}_{A} +
    \bs{v}^{\ast}_{B} + \bs{v}^{\ast}_{C}) \right)_1\right] + n \, y^{}_2
\left[\frac{b}{\ell} - \left( \pi_{\textnormal{int}}^{}(\bs{v}^{\ast}_{A} +
    \bs{v}^{\ast}_{B} + \bs{v}^{\ast}_{C}) \right)_2 \right] = 0  
\end{multline} 
because each of the terms in square brackets vanishes. Therefore
\begin{equation*}
c^{}_{k}(\varLambda^{}_{\varphi}(\Omega^{}_{i})) \, = \,
c^{}_{k'}(\varLambda^{}_{\varphi}(\Omega^{}_{i}))  
\end{equation*}
holds, and the spectrum is periodic with period $1/\ell$.

To obtain the diffraction spectrum of Kol$(3,1)$ from here, one only has to
rescale the positions in $\varLambda^{}_{\varphi}(\Omega)$ by a factor of
$1/\ell$. To summarize:

\begin{theorem}
The bi-infinite sequence Kol$(3,1)$, represented with equal bond lengths, is a
deformed model set and has a pure point diffraction spectrum.\qed  
\end{theorem}

\smallskip
\noindent
\textsc{Remarks}: By the same method, we can also find a deformation 
$\tilde{\varphi}(x^{\star}) = \tilde{a} \,x^{\star}_{1} + \tilde{b} \,
x^{\star}_{2}$ such that we represent the letter `$1$' of Kol$(3,1)$ with an
interval of length $\tilde{\ell}$ and the letter `$3$' with one of length $3 \,
\tilde{\ell}$. For this, the letters $A, B, C$ have bond lengths $6
\,\tilde{\ell},\,4\,\tilde{\ell},\,2\,\tilde{\ell}$, respectively. For the
parameters of the deformation (the average bond length must be $\ell$ again),
we get   
\begin{equation*}
\begin{split}
& \tilde{\ell} = \frac14 \, (7\, \alpha^2 - 15\,\alpha+1) \approx 0.49, \qquad
\tilde{a} = \frac12 \, (7\, \alpha^2 - 15\,\alpha-1) \approx -0.016 \quad
\text{and}\\ 
& \tilde{b} = \frac{\IM(\beta)}{59}\cdot(-179\,\alpha^2+379\,\alpha+3) =
\frac1{2\,\sqrt{59}^3} (1239\,\alpha^2-767\,\alpha-4661) \approx -0.36. 
\end{split}
\end{equation*}
Now (\ref{eq:invlen}) changes to
\begin{equation*}
\frac{n}{\tilde{\ell}} =
\pi(6 \, n \, \bs{v}^{\ast}_{A}+ 4 \, n \, \bs{v}^{\ast}_{B}+2 \, n \,
\bs{v}^{\ast}_{C}),  
\end{equation*}
and with the same calculation as before one gets an equation which corresponds
to (\ref{eq:bra2ket}), where the two terms in square brackets
\begin{equation*}
\left[\frac{\tilde{a}}{\tilde{\ell}} - \left(
    \pi_{\textnormal{int}}^{}(6\, \bs{v}^{\ast}_{A}+4\, \bs{v}^{\ast}_{B}+2\,
    \bs{v}^{\ast}_{C}) \right)_1\right] \quad \text{and} \quad
\left[\frac{\tilde{b}}{\tilde{\ell}} - \left(
    \pi_{\textnormal{int}}^{}(6 \,\bs{v}^{\ast}_{A}+4 \,\bs{v}^{\ast}_{B}+2
    \,\bs{v}^{\ast}_{C}) \right)_2 \right] 
\end{equation*}
also both vanish. Therefore, the spectrum is periodic with period
$1/\tilde{\ell}$ as expected~\cite{Baa01}, since
$\varLambda^{}_{\tilde{\varphi}}(\Omega) \subsetneqq \tilde{\ell} \cdot
\ZZ$. This representation with integer bond lengths (after rescaling) has the
advantage that the union of the three aperiodic sets
$\varLambda^{}_{\tilde{\varphi}}(\Omega^{}_{i})$ is still an aperiodic
set. Clearly, it is also pure point diffractive.

\smallskip
\noindent
Kol$(3,1)$ in its natural setting with intervals of length $1$, or of lengths
$3$ and $1$, can be obtained as a deformation of the model set $\CPS$ derived
above, where the intervals have incommensurate length. The basic theory of
this is fully developed in~\cite{BD00, Hof97}, but one can also understand,
from a dynamical systems point of view, which deformations are stable in the
sense that they do not change the spectral type of the dynamical spectrum (and
hence of the diffraction spectrum, due to unique ergodicity), see~\cite{BL03}.

\section*{Acknowledgments}

It is a pleasure to thank Christoph Bandt for fractal advice, Robert V.~Moody
for helpful discussions and the German Research Council (DFG) for financial
support. Also, we like to thank the referee for useful suggestions which led to
an improvement of this article.

\bigskip

\begin{appendix}
\section*{Appendix: An alternative proof of Theorem~\ref{thm:generic}}

By Proposition~\ref{prop:C0} we can choose a sequence $(c^{}_{i})^{}_{i}$ with 
$ c^{}_{i} = c+t^{\star}_{i} \in C^{}_0$ ($t^{\star}_{i} \in L^{\star}$) such
that $B^{}_{r^{}_{i}} \cap \CPS =  B^{}_{r^{}_{i}} \cap
\varLambda(c^{}_{i}+\Omega)$ for every sequence $(r^{}_{i})^{}_{i}$ with
$r^{}_{i} > i$. Also, this statement holds for every subsequence
$(c^{}_{i_{j}})^{}_{j}$. 
 
Now, assume $\partial\diamondset{\Omega} \cap L^{\star} \neq
\varnothing$. Then we have a point $\tilde{x} \in
\partial\diamondset{\Omega} \cap L^{\star}\cap\openset{\Omega}$ with
$\tilde{x} \in \partial\Omega^{}_{A} \cap \partial\Omega^{}_{B}$. Set
$\varepsilon^{}_0 = \operatorname{dist}(\tilde{x}, \partial\diamondset{\Omega}
\setminus (\partial\Omega^{}_{A} \cap \partial\Omega^{}_{B})) > 0$. Then a
translation $y \neq 0$ of $\Omega$ with $y \in B^{}_{\varepsilon^{}_0}(0)
\cap C$ has the following effect: $\tilde{x} \in y+\diamondset{\Omega}$,
because by the definition of $C$, $\tilde{x}$ cannot be on the boundary $y +
\partial\diamondset{\Omega}$, and by the choice of $\varepsilon^{}_0$, it
must either be in $y+\Omega^{}_{A}$ or $y+\Omega^{}_{B}$.  
    
Now take a sequence $(c^{}_{i})^{}_{i}$ as above. Clearly, this sequence must
converge to $0$. Therefore, there is an $N$ such that $|c^{}_{i}| <
\varepsilon^{}_0$ for all $i > N$. By choosing an appropriate subsequence
$(c^{}_{i_{j}})^{}_{j}$ we get a sequence $(\tilde{c}^{}_{j})^{}_{j}$ with
$\tilde{c}^{}_{j} = c^{}_{i_{j}}$ such that $\tilde{x}$ is always either in
$\tilde{c}^{}_{j}+\Omega^{}_{A}$ or in
$\tilde{c}^{}_{j}+\Omega^{}_{B}$. Also we have $B^{}_{r^{}_{i}} \cap \CPS =
B^{}_{r^{}_{i}} \cap \varLambda(\tilde{c}^{}_{i}+\Omega)$ for $r^{}_{i} >
i$. But both $\pi$ and $\pi^{}_{\textnormal{int}}$ are one-to-one. Therefore,
the inverse of the star map of $\tilde{x}$ must be a point of each
$\varLambda(\tilde{c}^{}_{i}+\Omega)$ and it also must be in $B^{}_{R}(0)$
for some $R < \infty$. But by Proposition~\ref{prop:inside}, it is not in
$\CPS$. Therefore we get a contradiction and our assumption is
wrong. So, $0 \in C$ and $\CPS$ is generic by
Proposition~\ref{prop:LI}\ref{en:C_generic}. \qed

\end{appendix}


\bigskip

\begin{thebibliography}{99}
\small

\bibitem{AI01}
P.~Arnoux and S.~Ito,
``Pisot substitutions and Rauzy fractals'',
\textit{Bull.\ Belg.\ Math.\ Soc.\ Simon Stevin} \textbf{8} (2001), 181--207.

\bibitem{Baa01}
M.~Baake,
``Diffraction of weighted lattice subsets'',
\textit{Canadian Math.\ Bulletin} \textbf{45} (2002), 483--498;
math.MG/0106111.

\bibitem{Baa02}
M.~Baake,
``A guide to mathematical quasicrystals'',
in: \textit{Quasicrystals}, eds.\ J.-B.~Suck,  M.~Schreiber and
P.~H\"{a}ussler, Springer, Berlin (2002), pp.\ 17--48; 
math-ph/9901014.

\bibitem{BL03}
M.~Baake and D.~Lenz,
``Deformation of Delone dynamical systems and topological conjugacy'';
in preparation.

\bibitem{BM00}
M.~Baake and R.V.~Moody,
``Self-similar measures for quasicrystals'',
in: \textit{Directions in Mathematical Quasicrystals},
eds.\ M.~Baake and R.V.~Moody,
AMS, Providence (2000), pp.\ 1--42; 
math.MG/0008063.

\bibitem{BM02}
M.~Baake and R.V.~Moody,
``Weighted Dirac combs with pure point diffraction'';
preprint 
math.MG/0203030.

\bibitem{BMS98}
M.~Baake, R.V.~Moody and M.~Schlottmann,
``Limit-(quasi)periodic point sets as quasicrystals with $p$-adic
  internal spaces'',
\textit{J.\ Phys.\ A: Math.\ Gen.} \textbf{31} (1998), 5755--5765; 
math-ph/9901008.

\bibitem{Ban97}
C.~Bandt,
``Self-similar tilings and patterns described by mappings'',
in: \textit{The Mathematics of Long-Range Aperiodic Order},
ed.\ R.V.~Moody, 
Kluwer, Dordrecht (1997), pp.\ 45--83.

\bibitem{BD00}
G.~Bernuau and M.~Duneau,
``Fourier Analysis of deformed model sets'',
in: \textit{Directions in Mathematical Quasicrystals},
eds.\ M.~Baake and R.V.~Moody,
AMS, Providence (2000), pp.\ 43--60.

\bibitem{BT86}
E.~Bombieri and J.E.~Taylor,
``Which distributions of matter diffract? An initial investigation'',
\textit{J.\ Physique Coll.}\ \textbf{C3} (1986), 19--29.

\bibitem{BS66}
S.I.~Borewicz and I.R.~\v{S}afarevi\v{c},
``Zahlentheorie'',
Birkh\"{a}user, Basel, 1966.

\bibitem{CS01}
V.~Canterini and A.~Siegel,
``Geometric representation of substitutions of Pisot type'',
\textit{Trans.\ Amer.\ Math.\ Soc.} \textbf{353} (2001), 5121--5144.

\bibitem{Dek78}
F.M.~Dekking,
``The spectrum of dynamical systems arising from substitutions of constant
length'',
\textit{Z.\ Wahrscheinlichkeitstheorie verw.\ Gebiete} \textbf{41} (1978),
221--239.  

\bibitem{Dek80}
F.M.~Dekking,
``Regularity and irregularity of sequences generated by automata'',
\textit{S\'{e}m.\ Th.\ Nombres Bordeaux} 1979--80, expos\'{e} 9, 901--910.

\bibitem{Dek97}
F.M.~Dekking,
``What is the long range order in the Kolakoski sequence?'',
in: \textit{The Mathematics of Long-Range Aperiodic Order},
ed.\ R.V.~Moody, 
Kluwer, Dordrecht (1997), pp.\ 115--125.

\bibitem{Edg90}
G.A.~Edgar,
``Measure, Topology and Fractal Geometry'',
Springer, New York, 1990.

\bibitem{GK97}
F.~G\"{a}hler and R.~Klitzing,
``The diffraction pattern of self-similar tilings'',
in: \textit{The Mathematics of Long-Range Aperiodic Order},
ed.\ R.V.~Moody, 
Kluwer, Dordrecht (1997), pp.\ 141--174.

\bibitem{Hof95}
A.~Hof,
``On diffraction by aperiodic structures'',
\textit{Commun.\ Math.\ Phys.}\ \textbf{169} (1995), 25--43.

\bibitem{Hof97}
A.~Hof,
``Diffraction by aperiodic structures'',
in: \textit{The Mathematics of Long-Range Aperiodic Order},
ed.\ R.V.~Moody, 
Kluwer, Dordrecht (1997), pp.\ 239--268.

\bibitem{Hut81}
J.E.~Hutchinson,
``Fractals and self-similarity'',
\textit{Indiana Univ.\ Math.\ J.}\ \textbf{30} (1981), 713--747.

\bibitem{Kol65}
W.~Kolakoski,
``Self generating runs, Problem 5304'',
\textit{Amer.\ Math.\ Monthly} \textbf{72} (1965), 674.

\bibitem{LM01}
J.-Y.~Lee and R.V.~Moody,
``Lattice substitution systems and model sets'',
\textit{Discrete Comput.\ Geom.}\ \textbf{25} (2001), 173--201;
math.MG/0002019.   

\bibitem{LMS}
J.-Y.~Lee, R.V.~Moody and B.~Solomyak,
``Pure point dynamical and diffraction spectra'',
\textit{Annales Henri Poincar\'e} \textbf{3} (2002), 1003--1018;
mp\_arc/02-39.

\bibitem{LGJJ93}
J.M.~Luck, C.~Godr\`{e}che, A.~Janner and T.~Janssen,
``The nature of the atomic surfaces of quasiperiodic self-similar
structures'',
\textit{J.\ Phys.\ A: Math.\ Gen.}\/ \textbf{26} (1993), 1951--1999.

\bibitem{Moo97}
R.V.~Moody,
``Meyer sets and their duals'',
in: \textit{The Mathematics of Long-Range Aperiodic Order},
ed.\ R.V.~Moody,
Kluwer, Dordrecht (1997), pp.\ 403--441.

\bibitem{Moo00}
R.V.~Moody,
``Model sets: a survey'', 
in: \textit{From Quasicrystals to More Complex Systems}, 
eds.\ F.~Axel, F.~D\'enoyer and J.P.~Gazeau,
EDP Sciences, Les Ulis, and
Springer, Berlin (2000), pp.\ 145--166;
math.MG/0002020.

\bibitem{Sch93}
M.~Schlottmann,
``Geometrische Eigenschaften quasiperiodischer Strukturen'',
Dissertation, Universit\"{a}t T\"{u}bingen (1993).

\bibitem{Sch98}
M.~Schlottmann,
``Cut-and-project sets in locally compact Abelian groups'',
in: \textit{Quasicrystals and Discrete Geometry},
ed.\ J.~Patera, 
AMS, Providence (1998), pp.\ 247--264.

\bibitem{Sch00}
M.~Schlottmann,
``Generalized model sets and dynamical systems'',
in: \textit{Directions in Mathematical Quasicrystals},
eds.\ M.~Baake and R.V.~Moody,
AMS, Providence (2000), pp.\ 43--60.

\bibitem{Sie02}
A.~Siegel,
``Repres\'{e}ntation des syst\`{e}mes dynamiques substitutifs non
unimodulaires'',
\textit{Ergodic Theory \& Dynam.\ Systems}, in press; preprint available 
from the author's homepage\footnote
{
Presently at:
  \texttt{http://www.irisa.fr/symbiose/people/siegel/Pro/publi.htm}
}.   

\bibitem{Diplom}
B.~Sing,
``Spektrale Eigenschaften der Kolakoski-Sequenzen'',
Diploma Thesis, Universit\"{a}t T\"{u}\-bingen (2002); 
available from the author.

\bibitem{InPrep}
B.~Sing,
``Kolakoski-$(2m,2n)$ are limit-periodic model sets'',
\textit{J.\ Math.\ Phys.} \textbf{44} (2003), 899-912;
math-ph/0207037.

\bibitem{Sir98}
V.F.~Sirvent,
``Mod\'{e}los geom\'{e}tricos asociados a substituciones'',
Habilitation (trabajo de ascenso), Universidad Sim\'{o}n Bol\'{i}var (1998);
available from the author's homepage\footnote
{
Presently at:
  \texttt{http://www.ma.usb.ve/\~{ }vsirvent/publi.html}
}.  

\bibitem{SS02}
V.F.~Sirvent and B.~Solomyak,
``Pure discrete spectrum for one-dimensional substitutions of Pisot type'',
\textit{Canadian Math.\ Bulletin} \textbf{45} (2002), 697--710. 

\bibitem{SW03}
V.F.~Sirvent and Y.~Wang,
``Self-affine tiling via substitution dynamical systems and Rauzy fractals'',
\textit{Pacific J.\ Math.} \textbf{206} (2002), 465--485.

\end{thebibliography}
\end{document}